\documentclass{siamart0516}

\usepackage{lineno}
\usepackage{amssymb,amsxtra,amsmath}
\usepackage{latexsym}
\usepackage{ifthen}
\usepackage{graphicx}
\usepackage{calc,fancyhdr}
\usepackage{indentfirst}
\usepackage{parskip,algorithm,algorithmicx}
\usepackage{algpseudocode}
\usepackage{color}
\usepackage{bigints}
\usepackage{subfig}
\usepackage{Article}
\usepackage{array,multirow,booktabs} % for fancier tables
\usepackage{mathtools}
\usepackage{cite}
\newcommand{\R}{\mathbb{R}}

\renewcommand{\B}{\mathbb{B}}

\newcommand{\revone}[1]{{\color{black} #1 }}
\newcommand{\revtwo}[1]{{\color{black}{#1}}}

\newcommand{\bs}[1]{\boldsymbol{#1}}

\newcommand{\TheTitle}{A Robust Hyperviscosity Formulation for Stable RBF-FD Discretizations of Advection-Diffusion-Reaction Equations on Manifolds} 
\newcommand{\TheAuthors}{Varun Shankar, Grady B. Wright, and Akil Narayan}

\headers{Hyperviscosity-based Stabilization for Transport on Manifolds}{\TheAuthors}

\title{{\TheTitle}\thanks{Submitted to the editors 08/01/2019.}}

\author{
  Varun Shankar\thanks{School of Computing, University of Utah, Salt Lake City, UT USA
    (\email{shankar@cs.utah.edu}). Corresponding author.}
  \and  
  Grady B. Wright\thanks{Department of Mathematics, Boise State University, ID, USA (\email{gradywright@boisestate.edu}).}
	\and
  Akil Narayan\thanks{Department of Mathematics, and Scientific Computing and Imaging (SCI) Institute, University of Utah, Salt Lake City, UT USA (\email{akil@sci.utah.edu}).}		
}

\usepackage[margin=1in]{geometry}
\hypersetup{
  pdftitle={\TheTitle},
  pdfauthor={\TheAuthors}
}
\begin{document}
\parindent=0pt
\maketitle

\begin{abstract}
  We present a new hyperviscosity formulation for stabilizing radial basis function-finite difference (RBF-FD) discretizations of advection-diffusion-reaction equations on manifolds $\mathbb{M} \subset \mathbb{R}^3$ of co-dimension one. Our technique involves automatic addition of artificial hyperviscosity to damp out spurious modes in the differentiation matrices corresponding to surface gradients, in the process overcoming a technical limitation of a recently-developed Euclidean formulation. Like the Euclidean formulation, the manifold formulation relies on von Neumann stability analysis performed on auxiliary differential operators that mimic the spurious solution growth induced by RBF-FD differentiation matrices. We demonstrate high-order convergence rates on problems involving surface advection and surface advection-diffusion. Finally, we demonstrate the applicability of our formulation to advection-diffusion-reaction equations on manifolds described purely as point clouds. Our surface discretizations use the recently-developed RBF-LOI method, and with the addition of hyperviscosity, are now empirically high-order accurate, stable, and free of stagnation errors.
\end{abstract}

\section{Introduction}
\label{sec:intro}

Radial Basis Functions (RBFs) are a powerful and flexible tool for generating numerical methods for the solution of partial differential equations (PDEs).  RBF collocation methods are easily applied to solving PDEs on irregular domains using scattered node layouts~\cite{Bayona2010,Davydov2011,Wright200699}. RBF-based methods also generalize naturally to the solution of PDEs on manifolds $\mathbb{M}\subset\mathbb{R}^3$ using only the Euclidean distance measure in the embedding space and Cartesian coordinates; see for example~\cite{FlyerWright:2007,FlyerWright:2009,FuselierWright:2013,Piret2012,FoL11,FlyerLehto2012,SWFKJSC2014,Piret2016,Aiton2011,LSWSISC2017, SWJCP2018,SNKJCP2018}. We focus exclusively on polynomially-augmented RBF-finite difference (RBF-FD) methods for their ability to harness the best features of both polynomial and RBF approximation.

It is natural to consider two types of (linear) stability in the context of RBF-FD methods: local stability, and global stability. Local stability refers to stability in the procedure of computing RBF-FD weights on each stencil. Global stability on the other hand refers to eigenvalue characterizations of differentiation matrices in RBF-FD discretizations of time-dependent PDEs: a globally-stable discretization for hyperbolic or elliptic operators would imply the absence of eigenvalues with positive real parts. Historically, RBF-FD methods have been locally unstable due to ill-conditioning in the RBF interpolation matrix~\cite{Fasshauer:2007,Wendland:2004}, manifesting as \emph{stagnation} in errors as the number of nodes is increased. Fortunately, work in the past two decades has helped overcome this on Euclidean domains, either by changing the RBF basis~\cite{FornbergPiret:2007,FaMC12,FLF,FoLePo13}, using rational approximations in the complex plane~\cite{FoWr,FoWr2016}, or by augmenting RBFs with high-degree polynomials~\cite{FlyerPHS,FlyerNS,FlyerElliptic}. On manifolds the local stability problem has been addressed either by using tangent plane projections on meshes to convert the problem to a Euclidean one~\cite{ReegerFornbergQuad1,ReegerFornbergQuad2,ReegerFornbergQuad3}, \revtwo{least-squares~\cite{Ling1}, the closest-point method which extends the problem to ambient space~\cite{Ling2}}, or a careful restriction of multivariate polynomials to the manifold~\cite{SNKJCP2018}. 

In the context of global stability, the augmenting of polyharmonic spline (PHS) RBFs with polynomials also appears to enhance global stability for elliptic operators on Euclidean domains~\cite{FlyerBoundary}, though this appears to depend on the node sets used (see~\cite{ShankarJCP2017,SFJCP2018} for evidence to the contrary of~\cite{FlyerBoundary}). However, the global stability problems in discretizations of hyperbolic operators are a topic of ongoing research. Unfortunately, while the RBF-LOI method presented in~\cite{SNKJCP2018} appears to be both locally and globally stable for elliptic operators on manifolds, it (like all RBF-FD methods) is not globally stable for hyperbolic operators. On the sphere (and in Euclidean domains), the current state-of-the-art is to add a stabilizing artificial hyperviscosity term of the form $\gamma_1 \Delta^{\gamma_2}$, \revone{$\gamma_1 \in \mathbb{R}$, $\gamma_2 \in \mathbb{N}$} an approach that was first used in the spectral methods literature under the title ``spectral (super) viscosity''~\cite{Tadmor89,MaSSV1,MaSSV2}. This was introduced to the RBF-FD context in~\cite{FoL11}, with empirical formulas for $\gamma_1$ and $\gamma_2$ provided in~\cite{FlyerLehto2012}. However, until recently, selecting $\gamma_1$ varied across applications~\cite{FlyerNS}.

In recent work~\cite{SFJCP2018}, the first author derived the first closed-form quasi-analytic formula for $\gamma_1$ and adapted a fully-analytic formula from~\cite{MaSSV1,MaSSV2} for $\gamma_2$ on Euclidean domains (applied to linear PDEs). These formulas generalize the spectral superviscosity formulas~\cite{Tadmor89,MaSSV1,MaSSV2} to scattered nodes and RBF-FD discretizations, and eliminates the need for parameter hand-tuning. The key contribution of that work was to use a 1D von Neumann analysis on a \emph{mathematical model} of the spurious growth modes of RBF-FD operators via so-called \emph{auxiliary derivatives} and \emph{auxiliary PDEs}. This approach allowed the author to derive quasi-analytic hyperviscosity formulas for stabilizing gradients and laplacians on Euclidean domains of arbitrary dimension. 

The goal of this article is to extend the Euclidean analysis of~\cite{SFJCP2018} to manifolds $\mathbb{M} \subset \mathbb{R}^3$, thereby developing a robust quasi-analytic hyperviscosity formulation for RBF-FD discretizations on manifolds. When combined with the recently developed RBF-LOI method, the resulting numerical method is free of stagnation errors, and is empirically both locally and globally stable for the examples we have investigated. To the best of our knowledge, this is the first globally-stable polynomially-augmented RBF-FD method \revtwo{of projection type} for hyperbolic problems on manifolds other than the sphere.  In addition, the use of the overlapped RBF-FD formulation~\cite{ShankarJCP2017,SFJCP2018,SNKJCP2018,SLMMS2018} makes the assembly of differentiation matrices efficient, especially for higher-order methods. \revone{Our hyperviscosity formulation is independent of RBF-LOI, and can be used straightforwardly with other RBF-FD methods based on the projection approach~\cite{SWFKJSC2014,SNKJCP2018}, the closest-point approach~\cite{Ling1,Ling2}, or the orthogonal gradients method~\cite{Piret2012, Piret2016}}.

The remainder of this paper is organized as follows. In the next section, we review the RBF-LOI method on manifolds. In Section 3 we review the Euclidean hyperviscosity formulation from~\cite{SFJCP2018}, and use it to motivate our new quasi-analytic hyperviscosity formulation for manifolds; we verify that our formulation is valid for important classes of multistage and multistep time integration schemes. We validate our methods on advection and advection-diffusion problems in Section 4 by measuring errors and convergence rates. In Section 5, we present applications of our method to solving nonlinear advection-diffusion-reaction equations on manifolds other than the sphere and torus. We conclude with a summary of our results and a discussion of future work in Section 6.

%\begin{table}[h!]
  %\begin{center}
  %\resizebox{\textwidth}{!}{
    %\renewcommand{\tabcolsep}{0.4cm}
    %\renewcommand{\arraystretch}{1.3}
    %{\scriptsize
    %\begin{tabular}{@{}cp{0.8\textwidth}@{}}
      %\toprule
      %$\bs{x}$ & Point in $\R^d$ \\
      %$X$ & Collection of $N$ points $\{\bs{x}_1, \ldots, \bs{x}_N \}$ in $\R^d$\\
      %$n$ & Local stencil size \\
      %$P_k$ & The collection of the $n$ nearest neighbors in $X$ of $\bs{x}_k \in X$ \\
      %$\mathcal{I}_j^k$ & The index in $X$ of the $j$th point in the set $P_k$, with $\mathcal{I}_1^k \equiv k$ \\
      %$\mathcal{G}^x$ & $x$-component of the surface gradient\\
      %$\rho_k$ & Width of $x_k$-centered stencil $P_k$\\
      %$\delta$ & Overlap parameter \\
      %$m$ & Polyharmonic spline degree \\
      %$R_k$ & Subset of $P_k$ dictated by overlap parameter $\delta$ \\
      %$V_P$ & Least orthogonal interpolant polynomial subspace associated to point set $P$ \\
      %$M$ & Number of polynomial augmentation terms used in RBF-FD procedure\\
      %$h^k_j$ & Ordered orthonormal basis elements for $V_{P_k}$ \\
			%$\tau$ & LOI tolerance parameter \\
      %$G_k^x$ & RBF-FD weights for the operator $\mathcal{G}^x$ on stencil $k$ \\
    %\bottomrule
    %\end{tabular}
  %}
    %\renewcommand{\arraystretch}{1}
    %\renewcommand{\tabcolsep}{12pt}
  %}
  %\end{center}
  %\caption{Notation used throughout this article.}\label{tab:notation}
%\end{table}

\section{RBF-LOI on Surfaces}
\label{sec:rbfreview}
This section explains our spatial discretization methodology, which is based upon three main ideas. 
\begin{itemize}
  \item Section \ref{sec:rbf-fd}: polynomially-augmented RBF-FD is used as the fundamental approximation~\cite{Wright200699,FlyerElliptic,SFJCP2018}.
  \item Section \ref{sec:loi}: the choice of basis for the polynomials is made through the least orthogonal interpolant (LOI) procedure to ensure polynomial unisolvency for nodes on manifolds \cite{SNKJCP2018}.
  \item Section \ref{sec:overlapped}: overlapping is used to reduce the number of local stencil approximations that must be computed~\cite{ShankarJCP2017,SFJCP2018,SNKJCP2018,SLMMS2018}.,
\end{itemize}
We refer to our algorithm, which is the combination of these ingredients, as RBF-LOI on surfaces.

\subsection{RBF-FD on surfaces}\label{sec:rbf-fd}
\label{sec:description}

Let $X = \{\vx_k\}_{k=1}^N$ be a set of nodes on $\mathbb{M}\subset\mathbb{R}^3$. Given a fixed $k \in \{1, 2, \ldots, N\}$, we will describe  RBF approximations to surface differential operators over a \textit{local stencil} of $n << N$ nodes for each $\vx_k$. To denote the stencil nodes we use index sets as follows: for each $\vx_k$, let $P_k$ consist of $\vx_k$ and its $n-1$ nearest neighbors (measured using the standard Euclidean distance), and let the set $\{\calI^k_1,\hdots,\calI^k_n\}$ denote the indices into the global node set $X$ of the nodes in $P_k$, with $\calI^k_1 = k$.  We refer to $P_k$ as the local stencil for $\vx_k$.  The standard RBF-FD approach constructs such a stencil for each $k = 1, \ldots, N$, while the overlapped RBF-FD approach can use the same stencil for more than one node in $X$, thus reducing the total number of stencils and substantially reducing the computational cost of constructing the RBF-FD differentiation matrices (see Section \ref{sec:overlapped}).  For simplicity in the remainder of this discussion, we work with stencil $k = 1$; hence all our quantities will feature sub/superscripts ``1" that can be replaced by ``$k$" to apply to a general stencil.

Suppose we wish to approximate the surface gradient $\nabla_{\mathbb{M}}$, defined in Cartesian coordinates as:
\begin{align}\label{eq:sgrad-op}
\nabla_{\mathbb{M}} = (I - \vn \vn^T)\nabla \eqqcolon [\calG^x,\calG^y,\calG^z]^T,
\end{align}
where $\vn$ is the \revone{unit} outward normal to the surface, and $\nabla$ is the $\mathbb{R}^3$ gradient. We will focus our discussion on the surface gradient; the surface Laplacian is then computed from the surface gradient in an iterated fashion as in~\cite{SWFKJSC2014,SNKJCP2018}.

Focusing on the $\calG^x$ component of $\nabla_{\mathbb{M}}$, the overlapped RBF-FD weights for every node in the stencil $P_1$ are computed using a \emph{family} of polynomially-augmented local RBF interpolants on $P_1$ parametrized by variables $\vx$ and $\vy$:
\begin{align}
s_1(\vx,\vy) = \sum\limits_{j=1}^n (w^x)_j(\vy) \|\vx - \vx_{\calI^1_j}\|^m + \sum\limits_{i=1}^{M} \lambda_i(\vy) h^1_i(\vx),
\label{eq:rbf_interp}
\end{align}
where all superscripts ``1" refer to the stencil index, $\| \cdot \|$ denotes the two norm, $m$ is odd and corresponds to the order of the PHS RBF, and $\left\{h^1_i(\cdot)\right\}_{i=1}^M$ is a basis for trivariate polynomials of degree $\ell$. In 3D Euclidean domains, a standard choice for the polynomial basis would be $M$ monomials corresponding to the total-degree $\ell$ of the trivariate polynomial subspace. On a (locally-) algebraic manifold, this choice typically results in a numerical lack of polynomial unisolvency on the local stencil. The RBF-LOI augmentation described in Section \ref{sec:loi} mitigates this issue. \revtwo{Regardless, the polynomial coefficients $\lambda_i(\vy)$ serve as Lagrange multipliers that enforce polynomial reproduction on the RBF-FD weights; this is discussed further below.}

Here $\vx$ refers to the nodes used to compute the weights, \revtwo{while $\vy$ refers to the location at which the weights $(w^x)_j(\vy)$ are computed}. In other words, each interpolant in the family is given by varying $\vy$, with the $n$ RBF-FD weights for that interpolant given by $(w^x)_j(\vy)$\footnote{One could also (equivalently) describe the procedure for computing $(w^x)_j(\vy)$ as interpolation of some function at the nodes $\vx$, following by evaluation at the nodes $\vy$~\cite{SWFKJSC2014}.}. In the standard RBF-FD method, $\vy \equiv \vx_1$ (the stencil center), which allows us to omit the $\vy$ parameter altogether; in the overlapped RBF-FD method $\vy$ allows for computation of weights in some radius around the stencil center (see Section \ref{sec:overlapped}). 

The overlapped RBF-FD weights for the operators $\calG^x$,$\calG^y$, and $\calG^z$ at all stencil points $\vx_{\calI^1_j}$ (every point $\vy \in P_1$) are computed by imposing certain conditions on an appropriate version of \eqref{eq:rbf_interp}; \emph{e.g.}, for $\calG^y$ we replace $(w^x)$ in \eqref{eq:rbf_interp} with $(w^y)$. Considering the operator $\calG^x$ without loss of generality, we impose the following two (sets of) conditions on the interpolants \eqref{eq:rbf_interp}:
\begin{subequations}
\begin{align}
s_1(\vx_{\calI^1_j},\vx_{\calI^1_i}) &=  \lf.\lf(\calG^x \|\vx - \vx_{\calI^1_j}\|^m\rt)\rt|_{\vx = \vx_{\calI^1_i}}, &i=1,\hdots,n, j=1,\hdots,n, \label{eq:interp_constraint}\\
\sum\limits_{j=1}^n (w^x)^1_j(\vx_{\calI^1_k}) h^1_i(\vx_{\calI^1_j}) &= \lf.\lf(\calG^x h^1_i(\vx)\rt)\rt|_{\vx = \vx_{\calI^1_k}}, & k=1,\hdots,n, i=1,\hdots,M. \label{eq:poly_constraint}
\end{align}
\end{subequations}
The first set of conditions enforces that $s_1(\vx,\vy)$ interpolate the derivatives of the PHS RBF at all the points in $P_1$. The second set of conditions enforces polynomial reproduction/exactness on \emph{all} the overlapped RBF-FD weights \revtwo{using the polynomial coefficients $\lambda_i(\vy)$}. If a degree-$\ell$ polynomial space is employed for $h^1(\vx)$, then $M = {\ell + d \choose d}$; we discuss choosing $\ell$ further below. The (family of) interpolants \eqref{eq:rbf_interp} and the two conditions \eqref{eq:interp_constraint}--\eqref{eq:poly_constraint} can be collected into the following block linear system:
\begin{align}
\begin{bmatrix}
A_1 & H_1 \\
H_1^T & O
\end{bmatrix}
\begin{bmatrix}
G^x_1 \\
\Lambda_1
\end{bmatrix}
=
\begin{bmatrix}
B_{A_1} \\
B_{H_1}
\end{bmatrix},
\label{eq:rbf_linsys}
\end{align}
where
\begin{align}
  (A_1)_{ij} &= \|\vx_{\calI^1_i} - \vx_{\calI^1_j} \|^m, & i,j&=1,\hdots,n, \\
  (H_1)_{ij} &= h^1_j(\vx_{\calI^1_i}), & i&=1,\hdots,n,\; j=1,\hdots,M,\\
  (B_{A_1})_{ij} &= \lf.\calG^x \|\vx - \vx_{\calI^1_j} \|^m \rt|_{\vx = \vx_{\calI^1_i}}, & i,j& =1,\hdots,n, \\
  \label{eq:h-diff} (B_{H_1})_{ij} &= \lf.\calG^x h^1_j(\vx)\rt|_{\vx = \vx_{\calI^1_j}}, & i&=1,\hdots,M,\; j=1,\hdots,n,\\
  O_{ij} &= 0, & i,j &= 1,\hdots,M.
\end{align}
$G^x_1$ is the $n \times n$ local matrix of overlapped RBF-FD weights for the operator $\calG^x$, with each column containing the RBF-FD weights for a point $\vy \in P_1$. Assuming $H_1$ has full column rank, then the linear system \eqref{eq:rbf_linsys} has a unique solution. The above procedure can be repeated with the operators $\calG^y$ and $\calG^z$ to obtain the local differentiation matrices $G^y_1$ and $G^z_1$. By construction, the \emph{columns} of $G^x_1$ and its counterparts populate the rows of the global differentiation matrices $G^x$, $G^y$, and $G^z$; for an assembly algorithm, see Algorithm 1 in~\cite{SFJCP2018}. As described above, the RBF-FD surface gradient weights are computed for every point in the stencil; This allows one to compute RBF-FD weights for the local surface Laplacian as 
\begin{align}
L_1 = \lf(G^x_1\tilde{G}^x_1\rt)^T + \lf(G^y_1\tilde{G}^y_1\rt)^T + \lf(G^z_1\tilde{G}^z_1\rt)^T,
\label{eq:surf_lap_loc}
\end{align}
where the $\tilde{G}^x_1$, $\tilde{G}^y_1$, and $\tilde{G}^z_1$ contain only the columns of $G^x_1$, $G^y_1$, and $G^z_1$ corresponding to the nodes at which weights are desired (some small neighborhood around the stencil center)~\cite{SNKJCP2018}.

The local stencil size $n$, the polyharmonic spline order $m$, the polynomial degree $\ell$, and the corresponding polynomial space dimension $M$ are not specified above. Our selection for these parameters is done mostly as in~\cite{SNKJCP2018,SFJCP2018}. We set $m=2\ell+1$ and set $\ell$ to $\ell= \xi + \theta - 1$, where $\theta$ is the order of the differential operator, and $\xi$ is the desired order of approximation. Once $\ell$ is computed, we use $M = {\ell + d \choose d}$, where $d$ is the dimension of the ambient space. Next, the stencil size is chosen to be $n = 2 M +1$ for advection problems~\cite{SNKJCP2018}, and $n = 2M + \lf\lfloor \log(2M)\rt\rfloor$ for mixed PDEs~\cite{SFJCP2018}, with $d$ the dimension of the ambient space (in this paper, $d=3$). 

The entire procedure above must be performed for each stencil; two practical issues that arise in such a procedure are that (a) the columns of $H_1$ may not be linearly independent and result in lack of uniqueness for the system \eqref{eq:rbf_linsys}, and (b) computing the weights proceeds by looping over each stencil, $k = 1, \ldots, N$, which can be quite expensive even if it must only be performed once (despite its formal $O(N)$ computational complexity). The next two sections make modifications to the above procedure that mitigate these two issues.

\subsection{The LOI-generated polynomial basis}\label{sec:loi}
Standard RBF-FD discretizations on manifolds can result in matrices $H_1$ in \eqref{eq:rbf_linsys} that are rank-deficient, making \eqref{eq:rbf_linsys} a non-invertible system. The cause of this issue is that if $X$ lies on a manifold, then local stencils often lie on an approximate algebraic variety, which destroys polynomial unisolvency on that stencil. 

The RBF-LOI procedure in \cite{SNKJCP2018} circumvents this issue by numerically constructing a tailored polynomial space on each stencil in a way that ensures unisolvency. The main idea is to start with the total polynomial space of degree $\ell$ described in the previous section and then perform adaptation on this space by identifying (numerical) rank deficiencies. These numerical rank deficiencies are eliminated to within a tunable tolerance $\tau$ by removing \revone{or} adding certain basis elements. We give precise choices for $\tau$ in section \ref{sec:results} detailing our numerical results and in Algorithm \ref{alg:hyp}.

We refer to \cite{SNKJCP2018,narayan_stochastic_2012} for a more complete algorithm description. In brief, the algorithm takes as input a local stencil (e.g., $P_1$) and an \textit{a priori} determined orthonormal polynomial basis, and outputs a polynomial space (along with new basis functions) ensuring unisolvency. The main workhorse is the least orthogonal interpolation (LOI) algorithm \cite{narayan_stochastic_2012}, which involves only standard numerical linear algebraic factorization operations and detects rank deficiency via the parameter $\tau$. The input orthonormal basis is chosen as the tensor-product Chebyshev polynomials on the Euclidean bounding box of stencil $P_1$. The LOI algorithm outputs the basis functions $\{h^1_j(\cdot)\}_{j=1}^M$, which are subsequently used in \eqref{eq:rbf_interp} and \eqref{eq:rbf_linsys}. We have observed that use of this procedure eliminates solvability issues for \eqref{eq:rbf_linsys}.

\subsection{Overlapped RBF-FD}\label{sec:overlapped}
Computing standard RBF-FD weights for each stencil $k$, $k = 1, \ldots, N$, can be very computationally expensive. An overlap strategy pioneered in \cite{ShankarJCP2017} substantially reduces this heavy computational cost.  We briefly review this strategy here.  As with Section \ref{sec:rbf-fd}, we specialize our discussion and our notation to stencil $k = 1$ and note that generalizations can be performed by replacing various sub/superscripts ``1" with ``$k$".
 
Let $\delta \in (0,1]$ denote the \emph{overlap parameter} and define the \emph{stencil retention distance} as 
\begin{align}
\rho_1 = (1-\delta)\max\limits_{j=1, \ldots, n} \|\vx_{\calI^1_1} - \vx_{\calI^1_j}\|
\end{align}
where $\|\cdot\|$ is the Euclidean norm in $\R^3$.  The parameter $\rho_1$ defines the radius of the ball $\mathbb{B}_1$ centered at the node $\vx_{\calI^1_1}$. Let $r_1$ denote the number of nodes in $P_1$ that lie in $\B_1$, and $R_1$ contain the set of \revtwo {global indices of these $r_1$ nodes from $P_1$}. Then the overlapped RBF-FD method involves computing RBF-FD weights for \textit{all} the nodes whose indices are in $R_1$.  Note that we use $R_1$ to determine the columns of $G_x^1$, $G_y^1$, and $G_z^1$ to use to form  $\tilde{G}^x_1$, $\tilde{G}^y_1$, and $\tilde{G}^z_1$ in \eqref{eq:surf_lap_loc}. 

The {overlapped} RBF-FD method makes one important modification to the standard RBF-FD method: to avoid computing multiple sets of RBF-FD weights for a node $\vx_k$, we also require that weights computed for some node $\vx_k$ never be recomputed by some other stencil $P_i, i\neq k$. Let $N_{\delta}$ denote the total number of stencils. For a quasi-uniform node set $X$, then $N_{\delta} = \frac{N}{p}$, where $p = \max \lf((1-\delta)^d n,1 \rt)$, and $d$ is the dimension (in the above discussion, $d = 3$). If $\delta=1$, this gives us $N_{\delta} = N$, recovering the standard RBF-FD method. However, if $\delta <1$, then $N_{\delta} \ll N$, giving a significant speedup over the standard RBF-FD method. For a more detailed complexity analysis, see~\cite{ShankarJCP2017}. Given the polynomial degree $\ell$ described in Section \ref{sec:rbf-fd}, we choose the overlap parameter $\delta$ as follows:
\begin{align}
 \delta =
  \begin{cases} 
      \hfill 0.7 \hfill & \text{ if $\ell \leq 4$}, \\
      \hfill 0.5 \hfill & \text{ if $4 < \ell \leq 6$}, \\
			\hfill 0.4 \hfill & \text{ if $\ell > 6$}. \\
  \end{cases}
\label{eq:overlap_ch}
\end{align}
We find that these values of $\delta$ result in stable differentiation matrices, while also facilitating the rapid assembly of these matrices.

\section{Hyperviscosity formulations on manifolds}
\label{sec:hyp_stab}
\begin{algorithm}
\caption{Hyperviscosity-stabilized RBF-FD Discretizations on Manifolds}
\label{alg:hyp}
\begin{algorithmic}[1]	
  \Statex{\bf Given}: $X = \{\vx_j\}_{j=1}^N$, a set of nodes on the manifold.
	\Statex{\bf Given}: $h$, the average separation distance between nodes.
	\Statex{\bf Given}: $\xi$, the desired order of approximation of the numerical method.
	\Statex{\bf Given}: $\nu$, the surface diffusion coefficient.
	\Statex{\bf Given}: $\vu$, a velocity field tangent to the surface at every point.
	\Statex{\bf Given}: $c_0(\vx) = c(\vx,0)$, an initial condition.
	\Statex{\bf Generate}: $\underline{C} \approx \lf.c(\vx,t)\rt|_X$, the numerical solution to $\frac{\partial c}{\partial t} + \nabla_{\mathbb{M}}\cdot(c\vu) = \nu \Delta_{\mathbb{M}} c$.
	\State Set the polynomial degree to $\ell = \xi + 1$ if $\nu \neq 0$, otherwise use $\ell = \xi$.
	\State Set the PHS RBF exponent to $m = 2\ell+1$.
	\State Let $M = {\ell + d \choose d}$, where $d$ is the spatial dimension.
	\State Set the stencil size to $n = 2M + \lf\lfloor \ln 2M \rt\rfloor$ if $\nu \neq 0$, otherwise set $n = 2M+1$.
	\State Set the hyperviscosity exponent to $\gamma_2 = \lf \lfloor \ln n\rt\rfloor$ if solution is smooth, otherwise set $\gamma_2 = 2$.
	\If {$\nu \neq 0$}	
	\State Set the LOI tolerance to $\tau = 10^{-3}$ if $\ell \leq 4$, else set $\tau = 10^{-4}$.
	\Else
	\State Set the LOI tolerance to $\tau = 0.05$ if $\ell < 4$, $\tau = 10^{-3}$ if $\ell \in [4,5]$, or $\tau = 10^{-4}$ if $\ell > 5$.
	\EndIf
	\State Set the overlap parameter $\delta$ according to \eqref{eq:overlap_ch}.
	\State Compute local differentiation matrices $G^x_j$, $G^y_j$, and $G^z_j$ using \eqref{eq:rbf_linsys} and its counterparts.
	\State Compute local differentiation matrices $L_j$ using \eqref{eq:surf_lap_loc}.
	\State Assemble local differentiation matrices into (sparse) global differentiation matrices $G^x$, $G^y$, $G^z$, and $L$ using Algorithm 1 from~\cite{SFJCP2018}.
	\State Compute (sparse) global differentiation matrix for hyperviscosity operator as $H = L^{\gamma_2}$.
	\State Estimate $\tau_x$, $\tau_y$, and $\tau_z$, the real parts of eigenvalues with largest real parts for $G^x$, $G^y$, and $G^z$, respectively, as described in Section \ref{sec:growth_exp}.
	\State Estimate growth exponents $q_x$, $q_y$, and $q_z$ using \eqref{eq:q_x} and its counterparts.
	\State Estimate $\eta(\vx)$ using \eqref{eq:eta}, then estimate $\bar{\eta}$ using \eqref{eq:eta_bar}.
	\State Evaluate components of $\vu$ ($u^x$, $u^y$, and $u^z$) on $X$ to obtain $\underline{U}^x$, $\underline{U}^y$, and $\underline{U}^z$.
	\State Compute hyperviscosity coefficient $\gamma_1$ using \eqref{eq:hyp_man_analytic_final} if $\nabla_{\mathbb{M}} \cdot \vu = 0$, else use \eqref{eq:hyp_man_analytic_final2}.	
	\State Solve 
	\begin{align}
	\frac{\partial \underline{C}}{\partial t} + G^x \lf(\underline{U}^x \circ \underline{C}\rt) + G^y \lf(\underline{U}^y \circ \underline{C}\rt) + G^z \lf(\underline{U}^z \circ \underline{C}\rt) = \nu L \underline{C} + \gamma_1 H \underline{C} 
	\end{align}
	using a suitable time integration scheme (typically an IMEX method~\cite{Ascher97}).
\end{algorithmic}
\end{algorithm}

This section describes our new hyperviscosity formulation for manifolds. First, we review our Euclidean formulation in Section \ref{sec:hyp_euc_stab}. Next, we present our new manifold hyperviscosity formulation for the surface advection equation (Section \ref{sec:hyp_stab_man_main}) and the surface advection-diffusion equation (Section \ref{sec:hyp_stab_man_main2}). We then discuss computing growth exponents for our model of spurious growth in Section \ref{sec:growth_exp}. Finally, we conclude by discussing both the selection of $\gamma_2$ in the $\gamma_1 \Delta_{\mathbb{M}}^{\gamma_2}$ hyperviscosity term (Section \ref{sec:k_sel}) and the numerical approximation of the operator itself (Section \ref{sec:hyp_disc}). Our methodology is summarized in Algorithm \ref{alg:hyp}, which uses the surface advection-diffusion equation as an example but also specifies choices for pure surface advection. In Algorithm \ref{alg:hyp}, an underlined quantity (such as $\underline{C}$ or $\underline{U}^x$) refers to an evaluation of the corresponding scalar field on a node set $X$, resulting in a vector (an array) of length $N$. In addition, Algorithm \ref{alg:hyp} uses the $\circ$ notation to denote the element-wise (or Hadamard) product between underlined quantities.

\subsection{Review of Euclidean hyperviscosity formulations}
\label{sec:hyp_euc_stab}

We now review the hyperviscosity formulation for Euclidean domains $\Omega \in \mathbb{R}^d$ developed in~\cite{SFJCP2018}. Let $c(x,t)$ be a scalar field satisfying the linear advection equation
\begin{align}
\frac{\partial c}{\partial t} + u \frac{\partial c}{\partial x} = 0.
\label{eq:adv_1d}
\end{align}
Let $c(x,t) = \hat{c}(t) e^{i \hat{k} x}$, where $\hat{k}$ is the wave number. In addition, let $c^{n+1} = c(t_{n+1},x)$ and $c^n = c(t_n,x)$. Using the traditional von Neumann analysis framework, we write $c^{n+1} = \varrho c^n$, where \revtwo{$\varrho$} is an amplification/growth factor, \revtwo{and} $|\varrho|\leq 1$ is necessary for time stability. Consider a Forward Euler\footnote{The reader may find it odd that the discussion uses Forward Euler, which is unstable on \eqref{eq:adv_1d} even in the absence of spurious growth modes. However, this was chosen purely for simplicity of analysis and to demonstrate how one would cancel the spurious growth modes. In Section \ref{sec:hyp_stab_man_main}, we extend the analysis to more commonly-used time integrators} \revtwo{discretization} of \eqref{eq:adv_1d}, with the above relations substituted in. This yields:
\begin{align}
\frac{\varrho - 1}{\triangle t} + u i \hat{k} = 0.
\label{eq:adv_1d_fe}
\end{align}
This equation of course shows that $|\varrho| >1$ for Forward Euler, indicating the need for stabilization or alternative time discretizations. However, the problem can actually be more severe after spatial discretization. In the case of RBF-FD discretizations, the \emph{numerical} gradient typically introduces a spurious growth mode into the solution. Consider the \emph{auxiliary derivative} (an analytical analogue to the numerical gradient) by its action on plane waves as
\begin{align*}
\frac{\tilde{\partial} c}{\partial x} = (i \hat{k} - \tau_x \hat{k}^{q_x}) e^{i \hat{k}x},
\end{align*}
which assumes that $u>0$. Here, $\tau_x \in \mathbb{R}$ can be thought of as the real part of the eigenvalue with the largest real part in the spectrum of the RBF-FD differentiation matrix. We focus on the case where $\tau_x > 0$, since $\tau_x \leq 0$ does not correspond to spurious growth. The quantity $q_x$ is a growth exponent that can be estimated analytically when $\tau_x$ is known. Consequently, in the fully discretized setting, \eqref{eq:adv_1d_fe} is transformed into:
\begin{align}
\frac{\varrho - 1}{\triangle t} + u i \hat{k} = u \tau_x \hat{k}^{q_x}.
\label{eq:adv_1d_fe_sp}
\end{align}
Our goal now is to eliminate the right hand side of \eqref{eq:adv_1d_fe_sp} to nullify the spurious growth mode. This can be done by adding a constant times the power of the Laplacian to the right hand side of \eqref{eq:adv_1d}, i.e.\ $\gamma_1 \Delta^{\gamma_2}$. Substituting in plane waves into the modified equation and dividing out the exponential gives the following equation for the amplification factor:
\begin{align}
\frac{\varrho - 1}{\triangle t} + u i \hat{k} = u \tau_x \hat{k}^{q_x} + \gamma_1 (-1)^{\gamma_2} \hat{k}^{2\gamma_2}.
\label{eq:adv_1d_fe_sp_hyp}
\end{align}
From this, one obtains a formula for $\gamma_1$:
\begin{align*}
\gamma_1 = (-1)^{1-\gamma_2} u \tau_x \hat{k}^{q_x-2\gamma_2}.
\end{align*}
On a give node set of spacing $h$, the largest wave number that can be represented is $\hat{k} = 2h^{-1}$. This substitution yields a simple formula for $\gamma_1$:
\begin{align*}
\gamma_1 = (-1)^{1-\gamma_2} u \tau_x 2^{q_x-2\gamma_2} h^{2\gamma_2-q_x},
\end{align*}

For Euclidean domains $\Omega \subset \mathbb{R}^3$, the advection equation is given by
\begin{align}
\frac{\partial c}{\partial t} + \vu \cdot \nabla c &= 0,
\label{eq:adv_3d}
\end{align}
assuming that $\nabla \cdot \vu = 0$. We now \revone{have}three differentiation matrices for the gradient ($G^x, G^y$, and $G^z$ respectively). Each of these has their own spurious growth factor $\tau_x$, $\tau_y$, and $\tau_z$ respectively, and the corresponding growth exponents $q_x$, $q_y$, and $q_z$. Taking into account these generalizations, we obtain a new formula for $\gamma_1$:
\begin{align}
\gamma_1 = (-1)^{1-\gamma_2} 3^{-\gamma_2} \lf(u_x\tau_x 2^{q_x-2\gamma_2} h^{2\gamma_2-q_x} + u_y\tau_y 2^{q_y-2\gamma_2} h^{2\gamma_2-q_y} + u_z\tau_z 2^{q_z-2\gamma_2} h^{2\gamma_2-q_z}\rt).
\end{align}
If one desires a precomputed/fixed formula for $\gamma_1$, it may be convenient to approximate the above formula by:
\begin{align}
\gamma_1 \approx (-1)^{1-\gamma_2} 3^{-\gamma_2} \|\vu\| \lf(\tau_x 2^{q_x-2\gamma_2} h^{2\gamma_2-q_x} + \tau_y 2^{q_y-2\gamma_2} h^{2\gamma_2-q_y} + \tau_z 2^{q_z-2\gamma_2} h^{2\gamma_2-q_z}\rt).
\label{eq:hyp_gamma_3d}
\end{align}
The above formula assumes that the node spacing $h$ is uniform along all the spatial directions, but this restriction can easily be removed. For quasi-uniform nodes on some domain $\Omega \subset \mathbb{R}^d$, we set $h = N^{-1/d}$, where $N$ is the number of nodes used for discretization. On the other hand, if the nodes are instead scattered, we set $h$ to the be smallest distance between all pairs of nodes. This formula is slightly more general than the formula given in~\cite{SFJCP2018}, but as outlined in that work, the exponents $q_x$, $q_y$, and $q_z$ can be estimated analytically on a given node once the spurious eigenvalues have been estimated. This analysis carries over to a wide class of explicit, implicit, and implicit-explicit time integrators. Unfortunately, it does not necessarily carry over to manifolds.

\subsection{Hyperviscosity for the surface advection equation}
\label{sec:hyp_stab_man_main}

We turn our attention to the \emph{surface} advection equation:
\begin{align}
\frac{\partial c}{\partial t} + \vu \cdot \nabla_{\mathbb{M}}c = 0,
\end{align}
where we assume for simplicity of analysis that $\nabla_{\mathbb{M}}\cdot \vu = 0$. We discuss the case of $\nabla_{\mathbb{M}}\cdot \vu \neq 0$ in Section \ref{sec:div_hyp}. Once again, the RBF-FD discretization of the surface gradient operator can be modeled by the following \emph{auxiliary PDE}:
\begin{align}
\frac{\partial c}{\partial t} + \vu \cdot \tilde{\nabla}_{\mathbb{M}}c = 0,
\label{eq:man_adv_sp}
\end{align}
where the \emph{auxiliary surface gradient} $\tilde{\nabla}_{\mathbb{M}}$ is defined by its action on $c(\vx,t)$ as
\begin{align*}
\tilde{\nabla}_{\mathbb{M}} &= [ \tilde{\calG^x}, \tilde{\calG^y}, \tilde{\calG^z}]^T,
\end{align*}
which can be written out component-wise as:
\begin{align}
\tilde{\calG^x} c &= (\calG^x - \tau_x \hat{k}_x^{q_x})c, &
\tilde{\calG^y} c &= (\calG^y - \tau_y \hat{k}_y^{q_y})c, &
\tilde{\calG^z} c &= (\calG^z - \tau_z \hat{k}_z^{q_z})c.
\label{eq:mod_surface_grad}
\end{align}
Here, $\hat{k}_x$, $\hat{k}_y$, and $\hat{k}_z$ are wavenumbers. Substituting this definition of the auxiliary surface gradient into \eqref{eq:man_adv_sp} gives
\begin{align}\label{eq:tau-def}
  \frac{\partial c}{\partial t} + \vu \cdot \nabla_{\mathbb{M}}c &= (\vu \cdot \boldsymbol{\tau})c, &
  \bs{\tau} &\coloneqq [\tau_x \hat{k}^{q_x}, \tau_y \hat{k}^{q_y}, \tau_z \hat{k}^{q_z}]^T.
\end{align}
Our approach to cancel out the spurious term on the right is to use artificial hyperviscosity based on the surface Laplacian $\Delta_{\mathbb{M}}$, since the RBF-FD discretization of this operator appears to be globally stable on quasi-uniform node sets (empirically speaking)~\cite{SNKJCP2018,SWFKJSC2014,LSWSISC2017}. This gives the following PDE:
\begin{align}
\frac{\partial c}{\partial t} + \vu \cdot \nabla_{\mathbb{M}}c = (\vu \cdot \boldsymbol{\tau})c + \gamma_1(\vx) \Delta_{\mathbb{M}}^{\gamma_2} c,
\label{eq:surf_adv_hyp}
\end{align}
where $\gamma_1(\vx)$ is some undetermined function which will eventually be reduced to a single constant $\gamma_1$.

\subsubsection{Explicit Runge-Kutta (RK) methods}
\label{sec:erk}

We will now present a hyperviscosity formulation applicable to all explicit $s-$stage RK methods. The growth factor $\varrho$ for an explicit $s$-stage Runge-Kutta (RK) \revone{method}takes the form of the degree-$s$ Taylor polynomial:
\begin{align*}
\varrho(z) = \sum\limits_{j=0}^s \frac{z^j}{j!}.
\label{eq:gen_rk_gf}
\end{align*}
The value of $z$ depends on the equation being solved. To see this, we will first examine the Forward Euler method applied to \eqref{eq:surf_adv_hyp}. First, we observe that the growth factor $\varrho(z)$ for forward Euler (RK1) is given by
\begin{align}
\varrho(z) = 1 + z.
\label{eq:fe_growth_fac}
\end{align}
We now substitute in a plane wave of the form $e^{i {\bf \hat{k}} \cdot \vx}$ into \eqref{eq:surf_adv_hyp}, where ${\bf \hat{k}} = [\hat{k}_x,\hat{k}_y,\hat{k}_z]^T$ is the vector of wavenumbers, and use the growth factor for forward Euler. This gives us an equation of the form:
\begin{align}
\frac{\varrho-1}{\triangle t} \hat{c} e^{i {\bf \hat{k}} \cdot \vx} + \vu^n \cdot \nabla_{\mathbb{M}} \hat{c} e^{i {\bf \hat{k}} \cdot \vx} = (\vu^n \cdot \boldsymbol{\tau}) \hat{c} e^{i {\bf \hat{k}} \cdot \vx} + \gamma_1(\vx) \Delta_{\mathbb{M}}^{\gamma_2} \hat{c} e^{i {\bf \hat{k}} \cdot \vx}.
\label{eq:fe_G}
\end{align}
To solve the above equation for $\varrho$, we first need to simplify the $\vu^n \cdot \nabla_{\mathbb{M}} \hat{c} e^{i {\bf \hat{k} \cdot \vx}}$ term. This can be written as:
\begin{align}
\vu^n \cdot \nabla_{\mathbb{M}} \hat{c} e^{i {\bf \hat{k} \cdot \vx}} = \vu^n \cdot i\begin{bmatrix}
(1-n_x^2) k_x - n_x n_y k_y - n_x n_z k_z \\
-n_y n_x k_x + (1 - n_y^2) k_y - n_y n_z k_z \\
-n_z n_x k_x - n_z n_y k_y + (1-n_z^2) k_z
\end{bmatrix} \hat{c} e^{i {\bf \hat{k} \cdot \vx}} = (\vu \cdot P i{\bf k}) \hat{c} e^{i {\bf \hat{k} \cdot \vx}}.
\label{eq:surf_grad_pw}
\end{align}
This simplification gives us
\begin{align}
\frac{\varrho-1}{\triangle t} \hat{c} e^{i {\bf \hat{k}} \cdot \vx} + (\vu^n \cdot P i{\bf k}) \hat{c} e^{i {\bf \hat{k} \cdot \vx}} = (\vu^n \cdot \boldsymbol{\tau}) \hat{c} e^{i {\bf \hat{k}} \cdot \vx} + \gamma_1(\vx) \Delta_{\mathbb{M}}^{\gamma_2} \hat{c} e^{i {\bf \hat{k}} \cdot \vx}.
\label{eq:fe_G_2}
\end{align}
In the Euclidean case $\mathbb{M} \equiv \mathbb{R}^d$, we could now divide out the plane wave $e^{i {\bf \hat{k}} \cdot \vx}$ to obtain an expression for $\gamma_1$ since the plane wave is an eigenfunction of the $\mathbb{R}^d$ Laplacian. Unfortunately, for submanifolds $\mathbb{M} \subset \mathbb{R}^d$, plane waves are \emph{not} eigenfunctions of the surface Laplacian $\Delta_{\mathbb{M}}$. If we were to instead use the $\mathbb{R}^3$ Laplacian $\Delta$ in place of $\Delta_{\mathbb{M}}$, we could simply use \eqref{eq:hyp_gamma_3d}. However, it may in general be difficult to compute a reasonable discretization of $\Delta$ when the only points available to us are on the manifold $\mathbb{M}$. While this has been attempted for non polynomially-augmented RBF-FD~\cite{FlyerLehto2012,FoL11}, our attempt at doing so using the RBF-LOI procedure in Section \ref{sec:rbfreview} resulted in an approximation to $\Delta$ on the manifold that had positive real eigenvalues.

Our proposed approach to solving for $\gamma_1(\vx)$ is to approximate the action of $\Delta^{\gamma_2}_{\mathbb{M}}$ on plane waves as a function of the action of $\Delta^{\gamma_2}$. More specifically, we write
\begin{align}
\Delta^{\gamma_2}_{\mathbb{M}} e^{i {\bf \hat{k}} \cdot \vx} \approx \eta(\vx) \Delta^{\gamma_2} e^{i {\bf \hat{k}} \cdot \vx} =
\eta(\vx)(-1)^{\gamma_2} \lf(\hat{k}^2_x + \hat{k}^2_y + \hat{k}^2_z\rt)^{\gamma_2} e^{i {\bf \hat{k}} \cdot \vx},
\label{eq:lb_to_l_map}
\end{align}
where $\eta(\vx)$ is some unknown function. This approach allows us to treat plane waves as eigenfunctions of the surface Laplacian, overcoming the primary hurdle in computing $\gamma_1(\vx)$.  We can now use \eqref{eq:lb_to_l_map} to solve for $\varrho$ in \eqref{eq:fe_G_2}, giving
\begin{align*}
\varrho \approx 1 - \triangle t (\vu^n \cdot P{i\bf \hat{k}}) + \triangle t (\vu^n \cdot \boldsymbol{\tau}) + \triangle t  \gamma_1 \eta(\vx) (-1)^{\gamma_2} \lf(\hat{k}^2_x + \hat{k}^2_y + \hat{k}^2_z\rt)^{\gamma_2}.
\end{align*}
Comparing the above to \eqref{eq:fe_growth_fac} gives us an expression for $z$:
\begin{align*}
z = - \triangle t (\vu^n \cdot P{i\bf \hat{k}}) + \triangle t (\vu^n \cdot \boldsymbol{\tau}) + \triangle t  \gamma_1 \eta(\vx)(-1)^{\gamma_2} \lf(\hat{k}^2_x + \hat{k}^2_y + \hat{k}^2_z\rt)^{\gamma_2}.
\end{align*}
For an $s$-stage explicit RK method, we therefore have:
\begin{align}
\varrho(z) = \sum\limits_{j=0}^s \frac{1}{j!}\lf(- \triangle t (\vu^n \cdot P{i\bf \hat{k}}) + \triangle t (\vu^n \cdot \boldsymbol{\tau}) + \triangle t  \gamma_1 \eta(\vx)(-1)^{\gamma_2} \lf(\hat{k}^2_x + \hat{k}^2_y + \hat{k}^2_z\rt)^{\gamma_2}\rt)^j.
\end{align}
We now split $z$ into two components, $z = z_1 + z_2$, each associated with different terms in the PDE:
\begin{align*}
  z_1 &\coloneqq -\triangle t (\vu^n \cdot P{i\bf \hat{k}}), & z_2 &\coloneqq \triangle t (\vu^n \cdot \boldsymbol{\tau}) + \triangle t  \gamma_1 \eta(\vx)(-1)^{\gamma_2} \lf(\hat{k}^2_x + \hat{k}^2_y + \hat{k}^2_z\rt)^{\gamma_2}
\end{align*}
The term $z_1$ is that obtained for the semi-discrete surface advection equation without growth or hyperviscosity, while $z_2$ captures the growth and hyperviscosity terms. In the absence of spurious growth and hyperviscosity, we have the following identity:
\begin{align}
\varrho(z) = \varrho(z_1) = \sum\limits_{j=0}^s \frac{z_1^j}{j!} = \frac{e^{z_1} \Gamma(s+1,z_1)}{s!},
\label{eq:G1}
\end{align}
where $\Gamma(a,b)$ is the upper incomplete gamma function. However, if we allow growth and hyperviscosity, we have
\begin{align}
\varrho(z) = \varrho(z_1,z_2) = \sum\limits_{j=0}^s \frac{(z_1 + z_2)^j}{j!} = \frac{e^{z_1 + z_2} \Gamma(s+1,z_1 + z_2)}{s!}.
\label{eq:G2}
\end{align}
Subtracting \eqref{eq:G1} from \eqref{eq:G2} defines a term $\tilde{\varrho}$:
\begin{align*}
\tilde{\varrho} \coloneqq \varrho(z_1,z_2) - \varrho(z_1) = \frac{e^{z_1}\lf(e^{z_2} \Gamma(s+1,z_1+z_2) - \Gamma(s+1,z_1) \rt)}{s!}.
\end{align*}
Clearly $\tilde{\varrho}$ is the contribution of growth and hyperviscosity to the growth factor $\varrho(z)$. If $\tilde{\varrho} \equiv 0$, we would recover the growth factor for the explicit $s$-stage RK method on the pure semi-discrete surface advection equation. $\tilde{\varrho}$ can be zeroed out by setting $z_2 = 0$, which in turn implies:
\begin{align*}
&\triangle t (\vu^n \cdot \boldsymbol{\tau}) + \triangle t \gamma_1 \eta(\vx)(-1)^{\gamma_2} \lf(\hat{k}^2_x + \hat{k}^2_y + \hat{k}^2_z\rt)^{\gamma_2} = 0, \\
&\implies \gamma_1(\vx) = \frac{(-1)^{1-\gamma_2} \vu \cdot \bs{\tau}}{\eta(\vx) \lf(\hat{k}^2_x + \hat{k}^2_y + \hat{k}^2_z\rt)^{\gamma_2} }.
\end{align*}
Using the definition of $\bs{\tau}$ in \eqref{eq:tau-def} and making the simplifying assumption $\hat{k} = \hat{k}_x = \hat{k}_y = \hat{k}_z = 2h^{-1}$ yields
\begin{align*}
\gamma_1(\vx) = \frac{(-1)^{1-\gamma_2}3^{-\gamma_2}}{\eta(\vx)} \lf(u_x \tau_x 2^{q_x-2\gamma_2} h^{2\gamma_2-q_x} + u_y  \tau_y 2^{q_y-2\gamma_2} h^{2\gamma_2-q_y}+ u_z \tau_z 2^{q_z-2\gamma_2} h^{2\gamma_2-q_z}\rt).
\end{align*}
To allow for precomputation, we once again use the 2-norm instead:
\begin{align*}
\gamma_1(\vx) \approx \frac{(-1)^{1-\gamma_2}3^{-\gamma_2} \|\vu\|}{\eta(\vx)} \lf(\tau_x 2^{q_x-2\gamma_2} h^{2\gamma_2-q_x} + \tau_y 2^{q_y-2\gamma_2} h^{2\gamma_2-q_y} + \tau_z 2^{q_z-2\gamma_2} h^{2\gamma_2-q_z}\rt).
\end{align*}
The function $\eta(\vx)$ is in general not known exactly, but it can be approximated on a given node set $X$ on the manifold $\mathbb{M}$ as follows:  Let $L$ be the approximation to the surface Laplacian computed as in Section \ref{sec:rbfreview}. Then we can write:
\begin{align}
  \eta(\vx_j) &= \frac{\lf(L^{\gamma_2} \lf.e^{i {\bf \hat{k}} \cdot \vx}\rt|_{X}\rt)_j}{\lf((-1)^{\gamma_2} \lf(\hat{k}^2_x + \hat{k}^2_y + \hat{k}^2_z\rt)^{\gamma_2} \lf.e^{i {\bf \hat{k}} \cdot \vx}\rt|_{X}\rt)_j}, & j&=1,\ldots,N,
	\label{eq:eta}
\end{align}
where $N$ is the total number of nodes. Here, $\hat{k}_x$, $\hat{k}_y$, and $\hat{k}_z$ each are set to $2h^{-1}$. If the nodes are quasi-uniform, we set $h = N^{-1/2}$ for manifolds $\mathbb{M} \subset \mathbb{R}^d$ of co-dimension one. Finally, we define a single constant $\gamma_1$ as
\begin{align*}
\gamma_1 = \mathbb{E}_X\lf(\gamma_1(\vx) \rt),
\end{align*}
which is to be read as the ``real-valued expectation/mean of $\gamma_1(\vx)$ over the node set $X$''. This can in turn be written as:
\begin{align}
\gamma_1 = \frac{(-1)^{1-\gamma_2}3^{-\gamma_2} \|\vu\|}{\bar{\eta}}\lf(\tau_x 2^{q_x-2\gamma_2} h^{2\gamma_2-q_x} + \tau_y 2^{q_y-2\gamma_2} h^{2\gamma_2-q_y} + \tau_z 2^{q_z-2\gamma_2} h^{2\gamma_2-q_z}\rt),
\label{eq:hyp_man_analytic_final}
\end{align}
where $\bar{\eta}$ is designed to be real-valued and positive, given explicitly by:
\begin{align}
\bar{\eta} = \frac{1}{N}\sum\limits_{j=1}^N \lf( |\operatorname{Re}(\eta(\vx_j))|\rt).
\label{eq:eta_bar}
\end{align}
Though we have made a series of simplifying assumptions, we now have a simple and computable formula for $\gamma_1$ using \eqref{eq:hyp_man_analytic_final} (and \eqref{eq:eta_bar}). As we have shown, this formula is applicable to all $s$-stage explicit RK methods and any linear methods whose growth factor is expressible as the series \eqref{eq:gen_rk_gf}. 

\subsubsection{Explicit linear multistep methods}
\label{sec:elmm}
Explicit linear multistep methods (LMMs) are a popular class of methods for solving time-dependent problems. While their stability regions are not as large as the corresponding explicit RK methods, explicit LMMs are still competitive from the computational cost standpoint due to fewer function evaluations. We will now verify the hyperviscosity formulation in the context of these methods.

Before proceeding, we remark that \eqref{eq:hyp_man_analytic_final} (derived for RK methods) could be obtained if we had directly set the approximation for $\Delta^{\gamma_2}_{\mathbb{M}}$ to be:
\begin{align}
\Delta^{\gamma_2}_{\mathbb{M}}  e^{i {\bf \hat{k}} \cdot \vx} \approx \bar{\eta} (-1)^{\gamma_2} \lf(\hat{k}^2_x + \hat{k}^2_y + \hat{k}^2_z\rt)^{\gamma_2} e^{i {\bf \hat{k}} \cdot \vx}.
\label{eq:lb_to_l_map2}
\end{align}
It is also useful to define an analogous map for the surface Laplacian itself:
\begin{align}
\Delta_{\mathbb{M}}  e^{i {\bf \hat{k}} \cdot \vx} \approx -\bar{\omega} \lf(\hat{k}^2_x + \hat{k}^2_y + \hat{k}^2_z\rt) e^{i {\bf \hat{k}} \cdot \vx},
\label{eq:lb_to_l_map3}
\end{align}
where $\bar{\omega}$ is computed analogously to $\bar{\eta}$. We use both these new maps going forward. Unlike the explicit RK methods, deriving the growth factor for explicit LMMs is a slightly more elaborate procedure. Consider the generic ODE
\begin{align*}
\frac{\partial c}{\partial t} = \mu c.
\end{align*}
The $s$-step explicit LMM (the Adams-Bashforth method of order $s$) applied to this ODE gives:
\begin{align*}
c^{n+1} = c^n + \triangle t \sum\limits_{j=0}^{s-1} \alpha_j \mu c^{n-j},
\end{align*}
where $\alpha_j$ are some real-valued coefficients. Now, let $c^{n-s+1} = \hat{c} e^{i \hat{k} \cdot \vx}$, and $c^{n-s+2} = \varrho c^{n-s+1}$, where $\varrho$ is the growth factor. Substituting into the $s$-step Adams-Bashforth method, we obtain a polynomial equation for the growth factor $G$:
\begin{align}
\varrho^s - \varrho^{s-1}\lf(1 + \alpha_0 \triangle t\mu \rt) - \varrho^{s-2}\alpha_1\triangle t \mu - \varrho^{s-3}\alpha_2 \triangle t \mu - \ldots - \alpha_{s-1} \triangle t \mu = 0.
\label{eq:ab_s_growth_generic}
\end{align}
In the context of the surface advection equation with spurious growth and hyperviscosity, we have 
\begin{align*}
\triangle t \mu &= z_1 + z_2, \\
z_1 &= - \triangle t (\vu^n \cdot P{\bf k}),\\
z_2 &= \triangle t (\vu^n \cdot \boldsymbol{\tau}) + \triangle t  \gamma_1 \bar{\eta}(-1)^{\gamma_2} \lf(\hat{k}^2_x + \hat{k}^2_y + \hat{k}^2_z\rt)^{\gamma_2}.
\end{align*}
Substituting the above equations into \eqref{eq:ab_s_growth_generic} gives:
\begin{align*}
\varrho^s - \varrho^{s-1}\lf(1 + \alpha_0 z_1 + \alpha_0 z_2\rt) - \varrho^{s-2} \alpha_1 (z_1 + z_2) - \varrho^{s-3} \alpha_2 (z_1 + z_2) - \ldots - \alpha_{s-1} (z_1 + z_2) = 0.
\end{align*}
We can collect all the $z_1$ and $z_2$ terms separately to yield
\begin{align*}
\lf[\varrho^s - \varrho^{s-1}\lf(1 + \alpha_0 z_1\rt) - \varrho^{s-2} \alpha_1 z_1 - \ldots - \alpha_{s-1} z_1 \rt] - z_2\lf[\varrho^{s-1} \alpha_0 + \varrho^{s-2} \alpha_1 + \ldots +\alpha_{s-1} \rt] = 0.
\end{align*}
The first bracketed term above is the equation for the growth factor $\varrho$ on the pure semi-discrete surface advection equation, while the second bracketed term accounts for growth and hyperviscosity. To cancel out the latter, we either require the growth factor $\varrho$ to satisfy the polynomial $\varrho^{s-1} \alpha_0 + \varrho^{s-2} \alpha_1 + \ldots +\alpha_{s-1} = 0$ or we can simply set $z_2 = 0$. Doing the latter yields
\begin{align*}
\triangle t (\vu^n \cdot \boldsymbol{\tau}) + \triangle t \gamma_1 \bar{\eta}(-1)^{\gamma_2} \lf(\hat{k}^2_x + \hat{k}^2_y + \hat{k}^2_z\rt)^{\gamma_2} = 0,
\end{align*}
which results in the same expression for $\gamma_1$ as in the case of explicit RK methods, once again matching \eqref{eq:hyp_man_analytic_final}.

\subsection{Hyperviscosity for the surface advection-diffusion equation}
\label{sec:hyp_stab_man_main2}

We now consider the surface advection-diffusion equation defined as
\begin{align*}
\frac{\partial c}{\partial t} + \vu \cdot \nabla_{\mathbb{M}} c = \nu \Delta_{\mathbb{M}} c.
\end{align*}
This is a PDE of mixed character, with the ratio $\frac{\|\vu\|}{\nu}$ determining whether the transport of $c$ is dominated by advection or diffusion (over some characteristic length scale). In regimes where both advection and diffusion play an important role in the transport, numerical solution of the advection-diffusion equation is often done via an implicit-explicit (IMEX) method~\cite{Ascher97}. These methods allow us to advance the stiff diffusion term implicitly in time, while advancing the advection term explicitly. When considering the \emph{stabilized} auxiliary surface advection-diffusion equation,
\begin{align*}
\frac{\partial c}{\partial t} + \vu \cdot \tilde{\nabla}_{\mathbb{M}} c = \nu \Delta_{\mathbb{M}} c + \gamma_1 \Delta^{\gamma_2}_{\mathbb{M}} c,
\end{align*}
we must now decide whether to treat the hyperviscosity term explicitly or implicitly in time. If the problem contains enough diffusion to impact transport but insufficient diffusion to stabilize the auxiliary surface gradient, it may be reasonable to add a hyperviscosity term and treat it implicitly along with the stiff diffusion term (for efficiency and stability). While our formula for $\gamma_1$ is appropriate for the explicit treatment of the hyperviscosity term, we have yet to show whether this formula is sufficient for stability if the hyperviscosity term is treated implicitly.  We will focus on a popular IMEX time integrator: the semi-implicit backward difference formula of order 2 (SBDF2)~\cite{Ascher97}. First, consider a generic ODE of the form
\begin{align*}
\frac{d c}{dt} = \mu_1 c + \mu_2 c.
\end{align*}
If we decide to treat the $\mu_1$ term explicitly, and the $\mu_2$ term implicitly, the SBDF2 discretization of this ODE is:
\begin{align*}
\frac{3 c^{n+1} - 4c^n + c^{n-1}}{2\triangle t} = 2\mu_1 c^n - \mu_1 c^{n-1} + \mu_2 c^{n+1}.
\end{align*}
Substituting in a plane wave expression and using von Neumann analysis, we get a quadratic equation for the growth factor $\varrho$:
\begin{align}
\lf(1 - \frac{2}{3}\mu_2 \triangle t\rt) \varrho^2 - \frac{4}{3}\varrho\lf(1 + \triangle t \mu_1\rt) + \frac{1}{3} \lf(1 + 2\triangle t \mu_1 \rt) = 0.
\label{eq:sbdf2_generic}
\end{align}
In the context of the surface advection-diffusion equation, $\mu_1$ represents the action of the surface advection and spurious growth operators on plane waves, and can be written using \eqref{eq:surf_grad_pw} as
\begin{align}
\mu_1 &= - \vu \cdot P i {\bf \hat{k}} + \vu \cdot \boldsymbol{\tau},
\label{eq:mu1}
\end{align}
where we now assume for simplicity that the velocity $\vu$ is constant in time. $\mu_2$ represents the action of the surface Laplacian and surface hyperviscosity operators on plane waves. To derive an expression for $\mu_2$, we now also need the approximation \eqref{eq:lb_to_l_map3}, which allows us to write $\mu_2$ as:
\begin{align}
\mu_2 = -\nu \bar{\omega} \lf(\hat{k}^2_x + \hat{k}^2_y + \hat{k}^2_z\rt) + \gamma_1 \bar{\eta} (-1)^{\gamma_2} \lf(\hat{k}^2_x + \hat{k}^2_y + \hat{k}^2_z\rt)^{\gamma_2}.
\label{eq:mu2}
\end{align}
Substituting \eqref{eq:mu1} and \eqref{eq:mu2} into \eqref{eq:sbdf2_generic}, collecting the surface advection and diffusion terms into a term $t_1$, and collecting the spurious growth and hyperviscosity terms into the term $t_2$, we have:
\begin{align*}
t_1 &= \lf(1 + \frac{2}{3}\nu\triangle t \bar{\omega} \lf(\hat{k}^2_x + \hat{k}^2_y + \hat{k}^2_z\rt)\rt)\varrho^2 - \frac{4}{3}\varrho\lf(1 - \vu \cdot P i {\bf \hat{k}} \rt) + \frac{1}{3} \lf(1 - 2\vu \cdot P i {\bf \hat{k}} \rt), \\
t_2 &= -\frac{2}{3}\triangle t \varrho^2 \gamma_1\bar{\eta} (-1)^{\gamma_2} \lf(\hat{k}^2_x + \hat{k}^2_y + \hat{k}^2_z\rt)^{\gamma_2} + \frac{2}{3} \triangle t \vu \cdot \boldsymbol{\tau} (1-2\varrho), \\
t_1 &+ t_2 = 0.
\end{align*}
Here, $t_1$ is obtained when discretizing the surface advection-diffusion equation, while $t_2$ contains all additional terms. We thus require that $t_2 = 0$, \emph{i.e.},
\begin{align*}
-\varrho^2 \gamma_1\bar{\eta} (-1)^{\gamma_2} \lf(\hat{k}^2_x + \hat{k}^2_y + \hat{k}^2_z\rt)^{\gamma_2} + \vu \cdot \boldsymbol{\tau} (1-2\varrho) = 0.
\end{align*}
The roots of this quadratic equation are given by:
\begin{align*}
\varrho = -  \frac{\vu \cdot \boldsymbol{\tau}}{\gamma_1\bar{\eta} (-1)^{\gamma_2} \lf(\hat{k}^2_x + \hat{k}^2_y + \hat{k}^2_z\rt)^{\gamma_2}} \pm \frac{ \sqrt{\lf(\vu \cdot \boldsymbol{\tau}\rt) \lf(\vu \cdot \boldsymbol{\tau} + \gamma_1\bar{\eta} (-1)^{\gamma_2} \lf(\hat{k}^2_x + \hat{k}^2_y + \hat{k}^2_z\rt)^{\gamma_2}\rt)}}{\gamma_1\bar{\eta} (-1)^{\gamma_2} \lf(\hat{k}^2_x + \hat{k}^2_y + \hat{k}^2_z\rt)^{\gamma_2}}.
\end{align*}
For time stability, we require that $|\varrho|\leq 1$. In addition to the CFL condition, this can be achieved if
\begin{align*}
\vu \cdot \boldsymbol{\tau} + \gamma_1\bar{\eta} (-1)^{\gamma_2} \lf(\hat{k}^2_x + \hat{k}^2_y + \hat{k}^2_z\rt)^{\gamma_2} = 0.
\end{align*}
This results in the following expression for $\gamma_1$ after substituting in the definition of $\boldsymbol{\tau}$, using the wavenumber-$h$ relationship, and pulling out the magnitude of the velocity:
\begin{align*}
\gamma_1 = \frac{(-1)^{1-\gamma_2}3^{-\gamma_2} \|\vu\|}{\bar{\eta}}\lf(\tau_x 2^{q_x-2\gamma_2} h^{2\gamma_2-q_x} + \tau_y 2^{q_y-2\gamma_2} h^{2\gamma_2-q_y} + \tau_z 2^{q_z-2\gamma_2} h^{2\gamma_2-q_z}\rt).
\end{align*}
Once again, this is identical to \eqref{eq:hyp_man_analytic_final}. We have thus demonstrated that our formula is applicable for both explicit and implicit discretizations of the hyperviscosity operator. The above analysis carries over to higher-order SBDF methods, though it is more tedious and requires computing the roots of cubic and higher-degree polynomial equations in $\varrho$. 

\subsection{Additional considerations}
\subsubsection{On stabilization from model diffusion} In the setting of the surface advection-diffusion equation, it is reasonable to wonder how much model diffusion is required to cancel out the spurious growth modes in the discretized advection operator. Fortunately, our framework for computing $\gamma_1$ also allows us to estimate a lower-bound on the required amount of stabilizing model diffusion. Consider the auxiliary surface advection-diffusion equation given by:
\begin{align*}
\frac{\partial c}{\partial t} + \vu \cdot \tilde{\nabla}_{\mathbb{M}} c = \nu \Delta_{\mathbb{M}} c.
\end{align*}
The goal is to estimate a lower-bound on $\nu$ that allows us to cancel out spurious growth modes. Proceeding with the von Neumann analysis for SBDF2 as before, now using the map \eqref{eq:lb_to_l_map3}, and requiring that $|G|\leq 1$, we obtain:
\begin{align}
\nu = \frac{-\|\vu\|}{3\bar{\omega}} \lf(\tau_x 2^{q_x - 2} h^{2-q_x} + \tau_y 2^{q_y - 2} h^{2-q_y} + \tau_z 2^{q_z - 2} h^{2-q_z}\rt).
\end{align}
We illustrate the utility of this formula with an example. Assume that all quantities except $\nu$ and $\bar{\omega}$ are $O(1)$. Then, $\nu \approx \frac{h}{6\bar{\omega}}$. In this setting, for a coarse node set on the sphere ($N = 2562$ icosahedral nodes), experiments show that $\bar{\omega} \approx 0.03$, which implies that $\nu \gtrsim 5N^{-1/2}$ is the least amount of model diffusion required to cancel spurious growth modes. In other words, if $\nu \gtrsim 5h$ and $\|\vu\| = O(1)$, no hyperviscosity is required for stabilization. A careful analysis of the model diffusion using our technique may obviate the need for numerical hyperviscosity. However, in this work, we primarily focus on parameter regimes requiring hyperviscosity.

\subsubsection{On divergent velocity fields} 
\label{sec:div_hyp}
Our analysis thus far assumed that $\nabla_{\mathbb{M}} \cdot \vu = 0$. We now deal with the more general case of a divergent velocity field. Consider now the auxiliary divergent surface advection equation (with added hyperviscosity)
\begin{align}
\frac{\partial c}{\partial t} + \vu \cdot \tilde{\nabla}_{\mathbb{M}} c = -c \tilde{\nabla}_{\mathbb{M}} \cdot \vu + \gamma_1 \Delta^{\gamma_2}_\mathbb{M} c.
\end{align}
Once again using the forward Euler discretization for simplicity, standard von Neumann analysis for the variable $c$ (together with \eqref{eq:lb_to_l_map2}) gives us
\begin{align}
\frac{\varrho-1}{\triangle t} + \vu^n \cdot P i{\bf k}  = - \tilde{\nabla}_{\mathbb{M}} \cdot \vu^n +  \vu^n \cdot \boldsymbol{\tau} + \gamma_1 \bar{\eta}(-1)^{\gamma_2} \lf(\hat{k}^2_x + \hat{k}^2_y + \hat{k}^2_z\rt)^{\gamma_2}. 
\end{align}
The first term in the r.h.s can be expanded using the definition of the growth model, giving us
\begin{align}
\frac{\varrho-1}{\triangle t} + \vu^n \cdot P i{\bf k}  = -\nabla_{\mathbb{M}} \cdot \vu^n + 2\vu^n \cdot \boldsymbol{\tau} + \gamma_1 \bar{\eta}(-1)^{\gamma_2} \lf(\hat{k}^2_x + \hat{k}^2_y + \hat{k}^2_z\rt)^{\gamma_2}.
\end{align}
Note the factor of $2$ in front of the $\vu^n$ term. We now want $\gamma_1$ so that:
\begin{align}
2\vu^n \cdot \boldsymbol{\tau} + \gamma_1 \bar{\eta}(-1)^{\gamma_2} \lf(\hat{k}^2_x + \hat{k}^2_y + \hat{k}^2_z\rt)^{\gamma_2}= 0.
\end{align}
This in turn implies that $\gamma_1$ is given as
\begin{align*}
\gamma_1 = \frac{2(-1)^{1-\gamma_2}3^{-\gamma_2}}{\bar{\eta}}\lf(u_x\tau_x 2^{q_x-2\gamma_2} h^{2\gamma_2-q_x} + u_y\tau_y 2^{q_y-2\gamma_2} h^{2\gamma_2-q_y} + u_z\tau_z 2^{q_z-2\gamma_2} h^{2\gamma_2-q_z}\rt),
\end{align*}
which can be made amenable to precomputation by setting
\begin{align}
\gamma_1 = \frac{2(-1)^{1-\gamma_2}3^{-\gamma_2}\|\vu\|}{\bar{\eta}}\lf(\tau_x 2^{q_x-2\gamma_2} h^{2\gamma_2-q_x} + \tau_y 2^{q_y-2\gamma_2} h^{2\gamma_2-q_y} + \tau_z 2^{q_z-2\gamma_2} h^{2\gamma_2-q_z}\rt).
\label{eq:hyp_man_analytic_final2}
\end{align}
It is straightforward to show that this analysis applies to all explicit RK and multistep methods, and the SBDF2 method for the surface advection-diffusion equation. Though we focus on divergence-free fields in this article, our experiments have shown that the above formula works well to stabilize the surface advection equation when divergent velocity fields are used.

\subsection{Computing growth exponents}
\label{sec:growth_exp}

We now discuss how to compute the growth exponents $q_x$, $q_y$, and $q_z$. While this is similar to the Euclidean case from~\cite{SFJCP2018}, the presence of a manifold adds certain complications. Let $X = \{x_j\}_{j=1}^N$ be a set of nodes on the manifold where we wish to compute the discrete differentiation matrices $G^x$, $G^y$, and $G^z$. Consider the function $f(\vx) = e^{i {\bf \hat{k}} \cdot \vx}$. Its analytical surface gradient is given by:
\begin{align}
{\bf g} (\vx) = \nabla_{\mathbb{M}} f(\vx) = P i {\bf \hat{k}} f(\vx) = [g^x(\vx), g^y(\vx), g^z(\vx)]^T.
\end{align}
Let $\underline{f}$, $\underline{g}^x$, $\underline{g}^y$, and $\underline{g}^z$ be the evaluations of $f(\vx)$, $g^x(\vx)$, $g^y(\vx)$, and $g^z(\vx)$ on the node set $X$. The \emph{approximate} surface gradient of $f(\vx)$ on the node set $X$ is given component-wise by matrix-multiplication with the differentiation matrices:
\begin{align}
\underline{\tilde{g}}^x = G^x \underline{f}, \ \underline{\tilde{g}}^y = G^y \underline{f}, \ \underline{\tilde{g}}^z = G^z \underline{f}.
\end{align}
From our growth model, we know that $G^x$, $G^y$, and $G^z$ are represented by the auxiliary differential operators $\tilde{\calG^x} $, $\tilde{\calG^y}$, and $\tilde{\calG^z}$ respectively. Consequently, ignoring truncation errors, the growth model gives us:
\begin{align}
\underline{\tilde{g}}^x = \underline{g}^x - \tau_x \hat{k}_x^{q_x}\underline{f}, \
\underline{\tilde{g}}^y = \underline{g}^y - \tau_y \hat{k}_y^{q_y}\underline{f}, \
\underline{\tilde{g}}^z = \underline{g}^z - \tau_z \hat{k}_z^{q_z}\underline{f}.
\end{align}
Focusing momentarily on the $x$-component without loss of generality, this allows us to write
\begin{align*}
\|\underline{g}^x - \underline{\tilde{g}}^x \| = \tau_x \hat{k}^{q_x} \|\underline{f}\|,
\end{align*}
with similar expressions for the $y$ and $z$ components. Since $\tau_x \|\underline{f}\| \neq 0$, we divide through by this quantity, take natural logarithms, and rearrange to obtain:
\begin{align}
q_x = \frac{\ln\lf(\|\underline{g}^x - \underline{\tilde{g}}^x \|\rt) - \ln\lf(\tau_x\|\underline{f}\| \rt)}{\ln \hat{k}}, \label{eq:q_x}
\end{align}
with similar expressions for $q_y$ and $q_z$. The usual substitution of $\hat{k} = 2h^{-1}$ gives us specific values of $q_x$, $q_y$, and $q_z$ for a given node set once the $\tau_x$, $\tau_y$, and $\tau_z$ have been computed. As mentioned previously, these represent the real part of the eigenvalue with the largest real part in the spectrum of $G^x$, $G^y$, and $G^z$ respectively. Notice that if we over/underestimate these quantities, the formula \eqref{eq:q_x}  and its counterparts automatically under/overestimate $q_x$, $q_y$, and $q_z$. This makes the formula for $\gamma_1$ robust to loose tolerances used in estimating $\tau_x$, $\tau_y$, and $\tau_z$. In our experiments, we use the following iterative procedure to estimate these spurious eigenvalues:
\begin{enumerate}
\item Attempt to estimate $\tau_x$, $\tau_y$, and $\tau_z$ using an implicitly-restart Arnoldi method with a loose tolerance of $10^{-3}$, looking for the eigenvalue with the largest real part.
\item If no eigenvalues are found to this tolerance, double the tolerance until one can be found.
\end{enumerate}

\subsection{Selecting $\gamma_2$}
\label{sec:k_sel}
In the spectral methods literature, Ma recommends that $\gamma_2$ in $\gamma_1 \Delta^{\gamma_2}$ be chosen according to the relation $\gamma_2 \leq O( \ln N)$, where $N$ is the total number of points used in the discretization~\cite{MaSSV1,MaSSV2}. In~\cite{SFJCP2018}, our Euclidean formulation used a version of this scaling law so that $\gamma_2 = \lf \lfloor 1.5 \ln n\rt\rfloor$, where $n$ is the RBF-FD stencil size. 

On manifolds, we set $\gamma_2 = \lf \lfloor \ln n\rt\rfloor$ if the solution is smooth. This allows $\gamma_1 \Delta^{\gamma_2}$ to approach the spectral case as $n \to N$, and filters higher frequencies as the stencil size is increased.  If the solution is not smooth, we set $\gamma_2 = 2$ and do not scale with increasing $n$. In our experiments, we discovered undesirable oscillations when using $\gamma_2 >3$ for test cases with non-smooth solutions. The rationale is simple: if the solution $c$ is not sufficiently smooth, $\gamma_1 \Delta^{\gamma_2} c$ may not even be continuous if $\gamma_2$ is large. In general, if the smoothness of the solution is not known, it may be possible to use the native space norm as a smoothness indicator~\cite{iske2002} and set $\gamma_2$ accordingly. We leave this approach for future work.

\subsection{RBF approximation to $\gamma_1\Delta_{\mathbb{M}}^{\gamma_2}$}
\label{sec:hyp_disc}
Following the procedure outlined in Section \ref{sec:rbfreview}, we compute $\Delta_{\mathbb{M}}^{\gamma_2}$ by using RBF-LOI to obtain the local surface Laplacian $L_j$ on each stencil.  We then assemble the local $L_j$ matrices into a global sparse matrix $L$ for the surface Laplacian and compute the differentiation matrix for the hyperviscosity operator as $H = L^{\gamma_2}$. This approach proved successful since the spectrum of $L$ was always well-behaved in our experiments (no eigenvalues with positive real parts)\footnote{We also explored approximating the hyperviscosity operator locally on each stencil and then assembling the local operators, but this approach failed for larger $k$ due to spurious positive real eigenvalues in the spectrum of the resulting differentiation matrix.}. This is in contrast to our Euclidean formulation, where we approximated $\Delta^{\gamma_2}$ directly on each stencil~\cite{SFJCP2018}. The advantage of the new approach is that it obviates the need to separately approximate the $H$ matrix. At first glance, the observant reader may notice that this approximation to $H$ is potentially a very low-order one (possibly only first-order accurate if $\gamma_2$ is large enough). However, since $\gamma_1 = O(h^{2\gamma_2})$, the scaling with $\gamma_1$ makes $\gamma_1 H$ a reasonable approximation to $\gamma_1 \Delta_{\mathbb{M}}^{\gamma_2}$ in the sense that $\gamma_1 H$ possess similar spectral properties. We leave the tackling of potential spurious eigenvalues in $L$ (and therefore $H$) to future work, but note that we did not encounter this problem for any of the test cases in this article.

\revtwo{\section{Role of node sets and parameters}
\label{sec:params}

Though our hyperviscosity formulation is automatic, the actual magnitude of the parameter $\gamma_1$ on a given node set is dictated by the behavior of the function $\eta(\vx)$ and the largest wavenumber $\hat{k} \approx 2h^{-1}$. In the case of $\eta(\vx)$, we advocated computing a real-valued average $\bar{\eta}$; in the case of $\hat{k}$, we approximated as $2h^{-1}$, where $h$ was chosen as $\frac{1}{\sqrt{N}}$, an \emph{average} measure of node spacing.  We now explore the impact of these choices on stability and accuracy.

\revone{\subsection{Stability on non-uniform node sets}
\label{sec:node_sets}
\begin{figure}[ht!]
\includegraphics[scale=0.4]{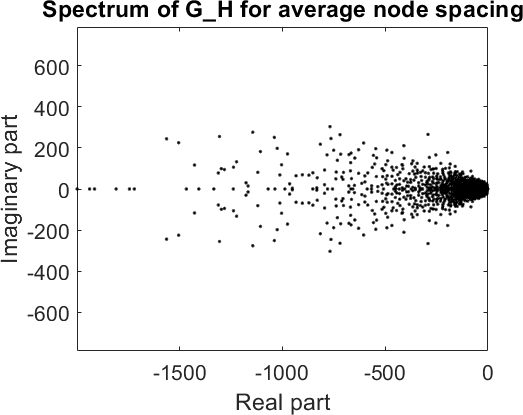}
\includegraphics[scale=0.4]{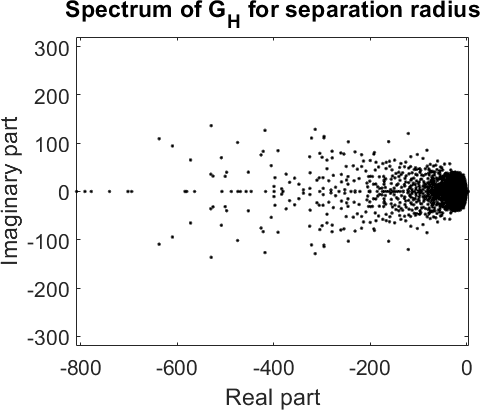}
\includegraphics[scale=0.4]{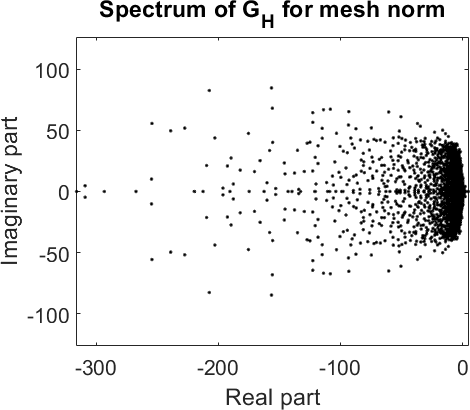}
\caption{Spectrum of $G_H$ for $N=1800$ non-uniform nodes on the torus when using $\bar{h}$ (left), $h_q$ (middle), and $h_{\rho}$ (right) in the formula for $\gamma_1$.}
\label{fig:eig_plots}
\end{figure}
We now explore the effect of these choices on the spectrum of the \emph{combined} discrete differentiation matrix $G_H = -G^x - G^y - G^z + \gamma_1 H$, assuming that the velocity is constant spatially. Ideally, the spectrum of this combined differential operator should have no eigenvalues with positive real parts. While the node sets used in Section 5 are all quasi-uniform, it is instructive to see the influence of non-uniformity on the value of $\gamma_1$ and hence the spectrum of $G_H$. First, we define three different quantities:
\begin{enumerate}
\item $\bar{h} = \frac{1}{\sqrt{N}}$, which we call the \emph{average node spacing}.
\item $h_q$, the \emph{separation radius}, defined as the minimal distance between any pair of points in the node set $X$\footnote{$h_q$ is referred to by the symbol $q$ in the RBF literature, and is often defined as half the value used in this work.}.
\item $h_{\rho}$, the \emph{mesh norm} or the \emph{fill distance}, defined as the radius of the largest ball that can be fit between any pair of nodes.
\end{enumerate}
In our formulas for $\gamma_1$, all three of these quantities are candidates for the variable $h$. For $N=2400$ (quasi-uniform) staggered nodes on a torus with inner radius $1/3$ and outer radius $1$, for instance, we have $\bar{h} = 0.0204$, $h_q = 0.0629$, and $h_{\rho} = 0.0871$. On the other hand, when we select a random subset of size $N=1800$ from these staggered nodes, the nodes are non uniform, and $h_{\rho}=0.1043$. Essentially, we expect $h_q$ and $h_{\rho}$ to be more different on non-uniform nodes. We set $\xi = 4$, and use Algorithm 1 to compute the different differentiation matrices. Finally, we compute the matrix $G_H$, and show its eigenvalues when each of the measures of node spacing are used above. The results are shown in Figure \ref{fig:eig_plots}. These results show that while all three measures of $h$ seem reasonable, the $h = \bar{h}$ (Figure \ref{fig:eig_plots} (left)) is the best choice if the goal is to guarantee the absence of an eigenvalue with any positive real part in the spectrum of $G_H$. Though it is hard to see, Figure \ref{fig:eig_plots} (middle) shows that using $h = h_q$ actually results in an eigenvalue with a very small real part, and  Figure \ref{fig:eig_plots} (right) shows that using $h = h_{\rho}$ results in an even larger real part. On a quasi-uniform node set, we have found that any of these measures of node spacing are appropriate, since they are all closer to each other. When solving PDEs, non-uniform node sets may arise in the context of $h$-refinement. Such non-uniform node sets are still more uniform than the node sets used in this experiment, so our experiment represents an extreme case. For the remainder of this article, we restrict our attention to quasi-uniform node sets and use $h = \bar{h}$.}

\subsection{Influence of $\eta(\vx)$ on stability}
\label{sec:eta}
\begin{figure}[ht!]
\includegraphics[scale=0.4]{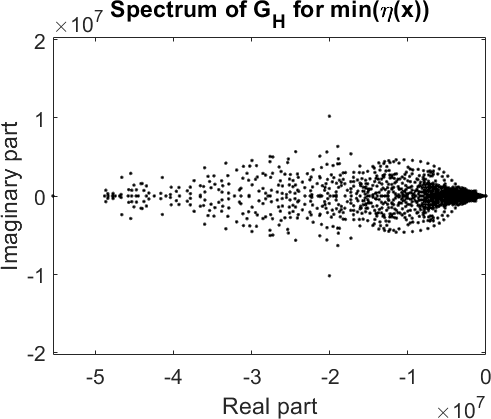}
\includegraphics[scale=0.4]{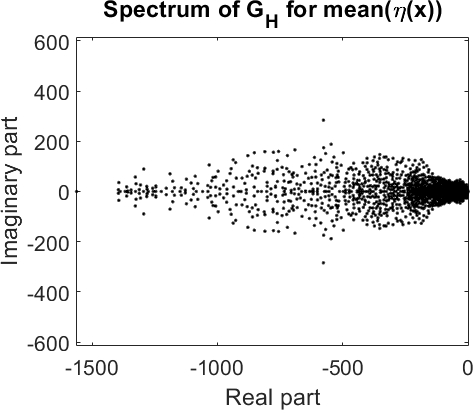}
\includegraphics[scale=0.4]{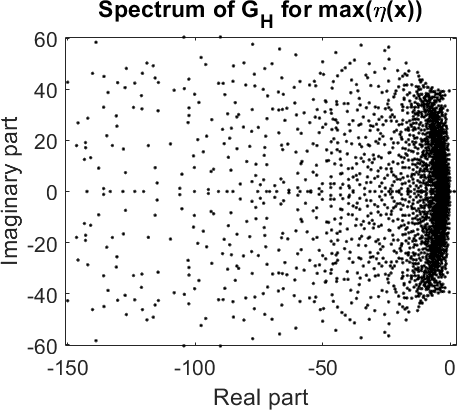}
\caption{Spectrum of $G_H$ for $N=2400$ quasi-uniform nodes on the torus when using $\min (\eta(\vx))$ (left), $\bar{\eta}$ (middle), and $\max (\eta(\vx))$ (right) in the formula for $\gamma_1$.}
\label{fig:eig_plots2}
\end{figure}
The formula for $\gamma_1$ contains the term $\bar{\eta}$, which is the real-valued expectation of the function $\eta(\vx)$. However, this choice was made simply to reduce $\eta(\vx)$ to a single positive number for computational convenience. It is possible to use other quantities, such as the maximum or minimum values. We now explore the influence of this choice on the spectrum of $G_H$. Once again, we set $\xi=4$ and use $N=2400$ quasi-uniform nodes on the torus. The results of the experiment are shown in Figure \ref{fig:eig_plots2}. The results are quite interesting. Using $\min (\eta(\vx))$ as in Figure \ref{fig:eig_plots2} (left) scatters the eigenvalues of $G_H$ very far into the left half of the complex plane, both along the real and imaginary axes (note the scale). While this differential operator is globally stable in the eigenvalue sense, our experiments indicate it is unusable for practical simulations as it significantly reduces the accuracy of our method. Figure \ref{fig:eig_plots2} (middle) shows that using $\bar{\eta}$ results in a much smaller scatter and has successfully shifted over eigenvalues with positive real parts. Finally, Figure \ref{fig:eig_plots2} (right) shows that using $\max (\eta(\vx))$ in fact produces the least scatter (again, note the scale), but produces an eigenvalue with a small positive real part (in this case, $\approx 1.5$). While $\bar{\eta}$ is clearly the safe choice, this indicates that $\max (\eta(\vx))$ can also be used in simulations if the time step $\Delta t$ is such that any small positive eigenvalue falls within the stability region of the time integrator.

\subsection{Effect of the analytic approximation to the eigenvalues of $\Delta_{\mathbb{M}}^{\gamma_2}$}
\label{sec:acc_eta}
Our technique for deriving an analytic expression for $\gamma_1$ relied on analytically approximating the action of $\Delta^{\gamma_2}_{\mathbb{M}}$ on plane waves as a function of the action of $\Delta^{\gamma_2}$. This was expressed through \eqref{eq:lb_to_l_map}, which we repeat here for convenience:
\begin{align*}
\Delta^{\gamma_2}_{\mathbb{M}} e^{i {\bf \hat{k}} \cdot \vx} \approx \eta(\vx) \Delta^{\gamma_2} e^{i {\bf \hat{k}} \cdot \vx} =
\eta(\vx)(-1)^{\gamma_2} \lf(\hat{k}^2_x + \hat{k}^2_y + \hat{k}^2_z\rt)^{\gamma_2} e^{i {\bf \hat{k}} \cdot \vx}.
\end{align*}
Here, $\eta(\vx)$ is computed numerically on a given node set $X$ by \eqref{eq:eta}. Our goal is to study the effect of this approximation when the eigenvalues of $\Delta^{\gamma_2}_{\mathbb{M}}$ are known. Let $f(\vx) = e^{i {\bf \hat{k}} \cdot \vx}$. To test this, we set $\gamma_2 = 1$, and compare $ \bar{\Delta}_{\mathbb{M}} f  = \bar{\eta} \lf(\hat{k}^2_x + \hat{k}^2_y + \hat{k}^2_z\rt) f$ with $L f$, where $L$ is the differentiation matrix that corresponds to ${\Delta}_{\mathbb{M}}$. Once again, we used staggered nodes on the torus for this experiment, with node sets of size $N=2400,5400,9600,21600,$ and $38400$, setting $\xi = 2,4,6$. We then measured the quantity $\|\lf(\bar{\Delta}_{\mathbb{M}} f\rt)_X - L f\|/\|L f\|$. We found that this was approximately constant ($\approx 0.44$) under both spatial refinement and order refinement with a magnitude that increased as $\gamma_2$ was increased. In fact, the quantity $\bar{\Delta}^{\gamma_2}_{\mathbb{M}} f$ appears to consistently overestimate $L^{\gamma_2} f = H f$, indicating that our formulation for $\gamma_1$ is likely overly strict, especially as $\gamma_2$ is increased. While this is reassuring for our time stability estimates for RK and IMEX methods, we plan in future work to explore spatially variable hyperviscosity formulations and analytic derivations for $\gamma_2$ to further automate these parameter choices. In addition, since the true surface Laplacian of $f$ can be written as $\Delta_{\mathbb{M}} f = \Delta f - 2H\frac{\partial f}{\partial {\vn}} - \frac{\partial^2 f}{\partial \vn^2}$, where $H$ is the mean curvature, we also plan to explore spatially-variable hyperviscosity formulations that take into account the curvature $H$.
}

\section{Results}
\label{sec:results}

\subsection{Surface Advection}

We now investigate the convergence properties of our stabilized RBF-FD methods on the surface advection equation. To the best of our knowledge, this is the first instance of a high-order RBF-FD method for transport on surfaces other than the sphere. On the sphere, we use two test cases: solid-body rotation of a cosine bell~\cite{Wil92} and deformational flow of two Gaussians~\cite{NairLauritzen2010}. On the torus, we measure the transport of both cosine bells and Gaussians in a time-independent, spatially-varying flow field. In all cases, the velocity fields return the initial conditions to their original spatial locations on the sphere and torus, allowing us to measure convergence rates. Our main focus is on convergence rates measured in the relative $\ell_2$-norm, measured by testing our methods on quasi-uniform node sets of increasing sizes. We use icosahedral nodes on the sphere, and staggered nodes on the torus, both as used in~\cite{SWFKJSC2014,LSWSISC2017,SNKJCP2018}. In all cases, we set $\xi$, the desired order of a accuracy for the spatial derivatives, equal to $\ell$ (since a first order differential operator is being approximated), and use $\xi = 2,4,6$. All subsequent parameters (including the overlap parameter $\delta$) were set using Algorithm \ref{alg:hyp}. For both tests, we use the classical third-order explicit Runge Kutta method (RK3). We chose RK3 instead of the more popular explicit RK4 to demonstrate that our hyperviscosity formulation does not rely on the larger stability region of explicit RK4. 

\subsubsection{Advection on the sphere}
\begin{figure}[h!]
\centering
\includegraphics[scale=0.6]{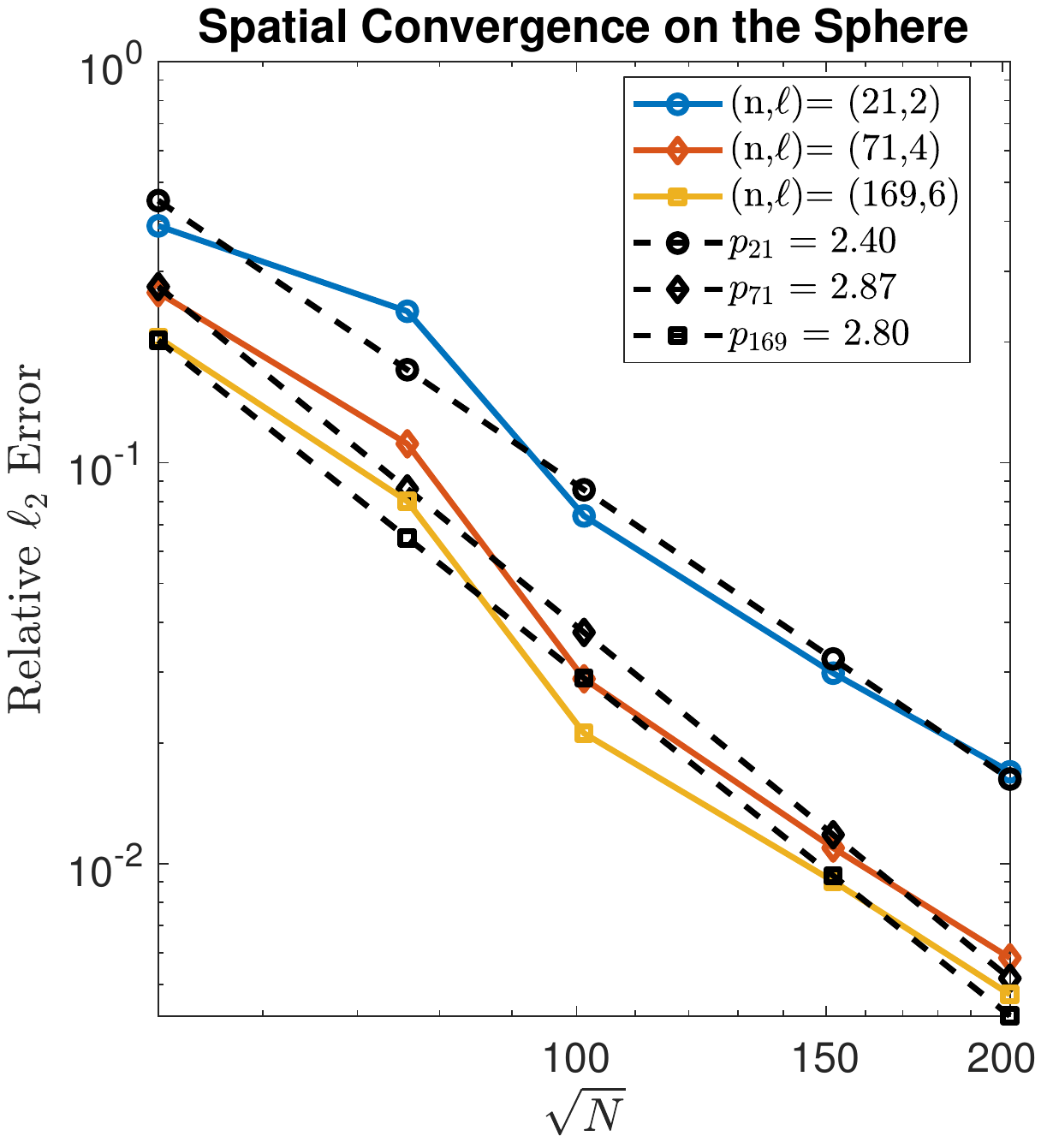}
\includegraphics[scale=0.6]{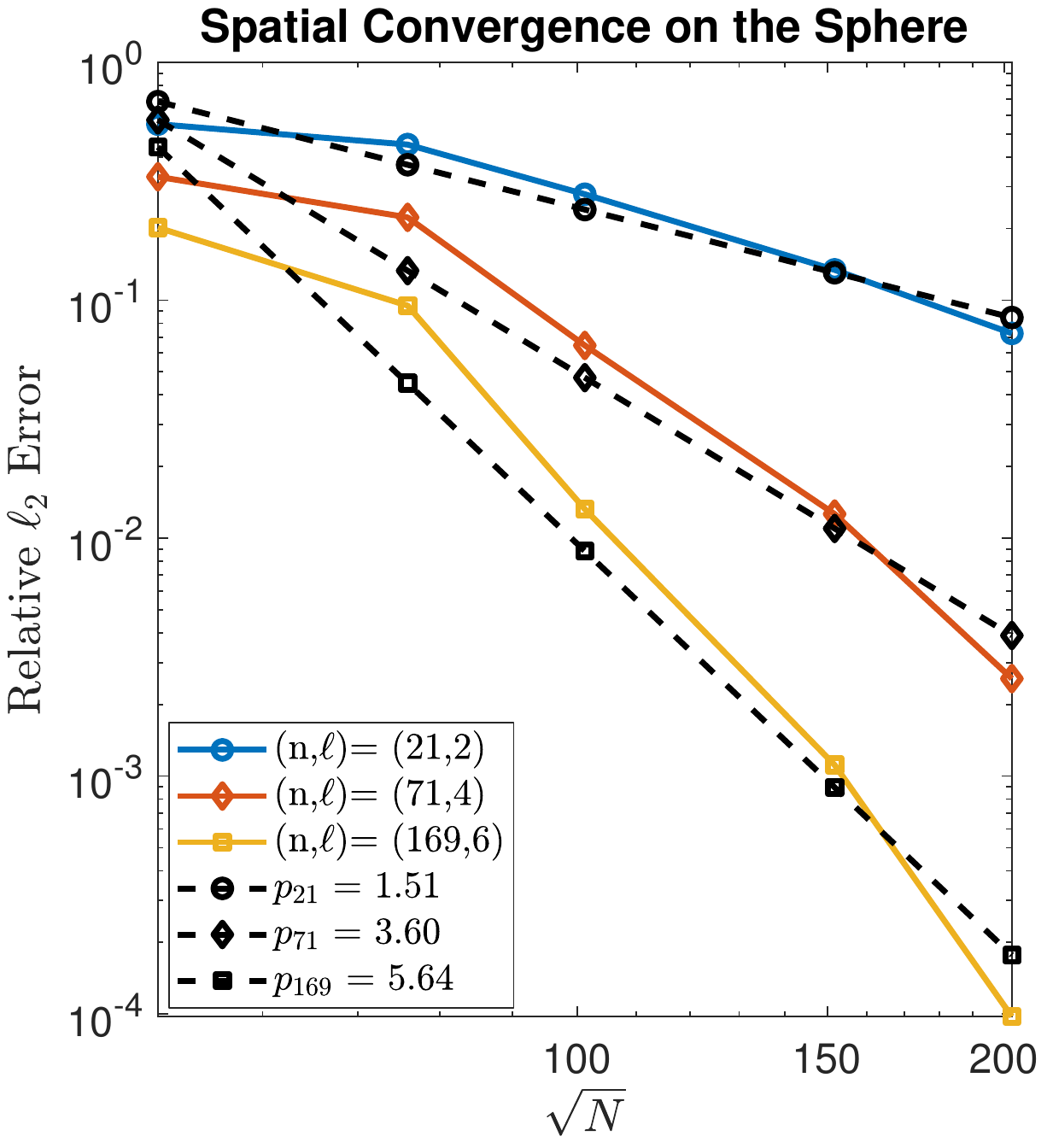}
\caption{Convergence on the sphere for the surface advection equation of a cosine bell in a steady flow (left) and gaussian in a deformational flow (right). The figure shows relative $\ell_2$-error as a function of $\sqrt{N}$ for different values of stencil $n$ and polynomial degree $\ell$. The dashed lines indicate lines of best fit, and their slopes are indicated in the legend as a measure of convergence rates.}
\label{fig:sphadv}
\end{figure}
We focus on two test cases for advection on the sphere. In both test cases, for a node set with $N$ nodes, we set the time step to $\triangle t = \frac{0.3}{\sqrt{N}}$, which roughly corresponds to a Courant number of $0.3$. This time step was chosen to approximately allow spatial errors to dominate over temporal errors. 

\paragraph{Solid body rotation of a cosine bell}
Our first test on the sphere is the solid body rotation test case from~\cite{Wil92}, which involves advection of a cosine bell in a steady velocity field. The components of the steady velocity field for this test case in spherical coordinates ($-\pi \leq \lambda \leq \pi$, $-\pi/2 \leq \theta \leq \pi/2$) are given by
\begin{align}
u(\lambda, \theta) =  \sin(\theta)\sin(\lambda)\sin(\alpha)-\cos(\theta)\cos(\alpha),\quad v(\lambda, \theta) =  \cos(\lambda)\sin(\alpha),
\label{eq:steady_vel}
\end{align}
where $\alpha$ is the angle of rotation with respect to the equator. We set $\alpha = \pi/2$ to advect the initial condition over the poles, and use a change of basis to obtain the velocity field in these coordinates.  The initial condition is a compactly-supported cosine bell centered at $(1,0,0)$:
\begin{align*}
c(\vx,0) = 
\begin{cases}
\frac{1}{2}\lf(1 + \cos\lf(\frac{\pi r}{R_{\rm b}} \rt) \rt) & \text{if}\ r < R_{\rm b}, \\
0 & \text{if}\ r \geq R_{\rm b},
\end{cases}
\end{align*}
where $\vx = (x,y,z)$, $r = \arccos(x)$, and $R_{\rm b}=1/3$. This initial condition is $C^1(\mathbb{S}^2)$, and one full revolution of the initial condition over the sphere requires simulation up to time $T = 2\pi$. The results of this simulation are shown in Figure \ref{fig:sphadv} (left), which shows that increasing the polynomial degree $\ell$ does not increase convergence rates when advecting the cosine bell, simply due to the limited smoothness of the cosine bell. However, as has been observed in the literature~\cite{SWJCP2018,FlyerLehto2012}, increasing $\ell$ does improve the accuracy.

\paragraph{Deformational flow of Gaussian bells}
The second test case involves advecting two Gaussian bells defined by
\begin{align*}
c(\vx,0) = 0.95\lf(e^{-5\|\vx - {\bf p}_1\|_2^2} + e^{-5\|\vx - {\bf p}_2\|_2^2} \rt),
\end{align*}
in a time-dependent deformational velocity field~\cite{NairLauritzen2010} whose components in spherical components are:
\begin{align*}
u(\lambda,\theta,t) &= \frac{10}{T} \cos\lf(\frac{\pi t}{T}\rt) \sin^2\lf(\lambda - \frac{2\pi t}{T} \rt)\sin\lf(2 \theta\rt) + \frac{2\pi}{T} \cos\lf(\theta\rt), \\
v(\lambda,\theta,t) &= \frac{10}{T} \cos\lf(\frac{\pi t}{T}\rt) \sin\lf(2\lambda - \frac{2\pi t}{T} \rt)\cos\lf(\theta\rt).
\end{align*}
Here, ${\bf p}_1 = \lf(\sqrt{3}/2,1/2,0\rt)$ and ${\bf p}_2 = \lf(\sqrt{3}/2,-1/2,0\rt)$. The flow field returns the solution to its initial position at final time $T=5$. This solution is $C^{\infty}(\mathbb{S}^2)$, and is ideal for measuring convergence rates. The convergence results in the $\ell_2$-norm are shown in Figure \ref{fig:sphadv}, which shows that increasing $\ell$ gives close to predicted convergence rates when advecting the Gaussian bells, despite the complexity of the deformational flow test case.

\subsubsection{Advection on the torus}
\begin{figure}[h!]
\centering
\includegraphics[width=0.3\textwidth]{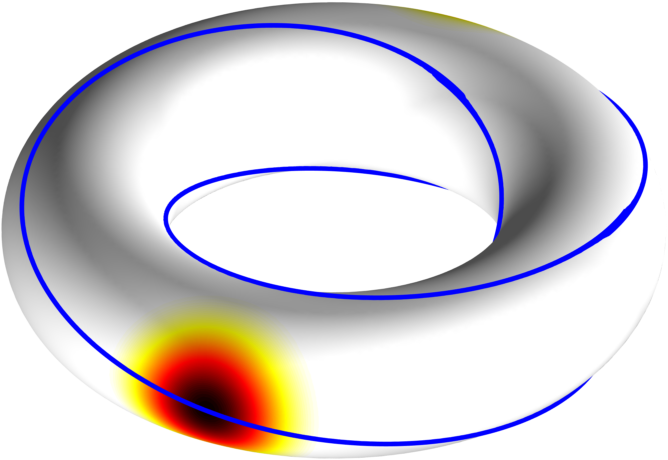}
\caption{Cosine bell initial condition for the torus test together with the path the solution advects over (according to \eqref{eq:sv_tor1}) superimposed on top in blue.}
\label{fig:toruscosinebell}
\end{figure}

For advection on the torus, we use the same velocity field for both tests. For a given set of points on the torus with inner radius $\frac{1}{3}$ and outer radius $1$ in Cartesian coordinates $(x,y,z)$, we define the following parametric coordinates:
\begin{align*}
&\rho = \sqrt{x^2 + y^2}, \ \phi = \frac{1}{3}\tan^{-1}\lf(\frac{y}{x}\rt), \ \theta = \phi - \frac{1}{2} \tan^{-1}\lf(\frac{z}{1 - r}\rt),\\
& \ \rho_1 =  1 + \frac{1}{3}\cos\lf(2\lf( \phi - \theta \rt) \rt),\ \rho_2 = -2r_1 \sin \lf(2 \lf( \phi - \theta \rt)\rt).
\end{align*}
Using these value, we then define the following time-independent tangential velocity field $\vu = (u,v,w)^T$ on the torus with respect to the Cartesian basis:
\begin{equation}
\begin{aligned}
u(\phi,\theta) &= \rho_1 \cos\lf(3\phi \rt) - \rho_2\lf(3\sin\lf(3\phi\rt) \rt),\\
v(\phi,\theta) &= \rho_1 \sin\lf(3\phi \rt) + \rho_2\lf(3\cos\lf(3\phi\rt) \rt),\\
w(\phi,\theta) &= -\frac{2}{3}\cos\lf(2 \lf(\phi -\theta \rt)\rt).
\end{aligned}
\label{eq:sv_tor1}
\end{equation}
For a particle placed at the initial position $(x_0,y_0,z_0)$ on the torus, this velocity field will advect the particle around a $(3,2)$ torus knot, returning the particle back to its initial position at $T=2\pi$ time units.  See Figure \ref{fig:toruscosinebell} for an illustration of the solution path.

We advect two different bell-type initial conditions, similar to the sphere, in this velocity field. The centers of these initial conditions are at the locations  $\vp_1 = \lf(1+1/3,0,0\rt)$ and $\vp_2 = -\vp_1$. In both test cases, for a node set with $N$ nodes, we set the time step to $\triangle t = \frac{0.3}{u_{max}\sqrt{N}}$, which roughly corresponds to a Courant number of $0.3$; here, $u_{max} = 4.1$ is the maximum pointwise magnitude of the velocity field over space. Again, this time step was chosen to approximately allow spatial errors to dominate over temporal errors introduced by RK3 integration. 

\begin{figure}[h!]
\centering
\includegraphics[scale=0.6]{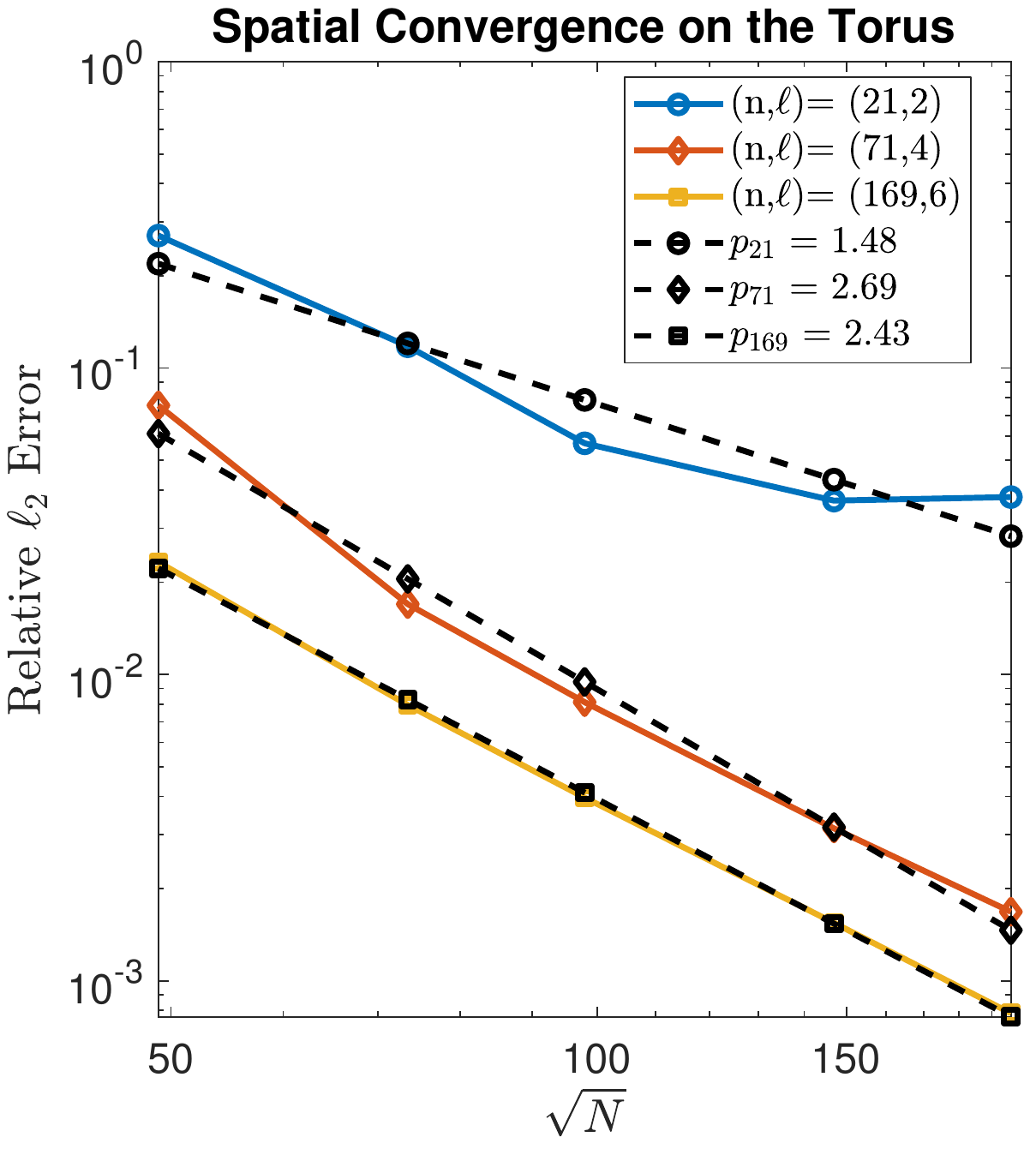}
\includegraphics[scale=0.6]{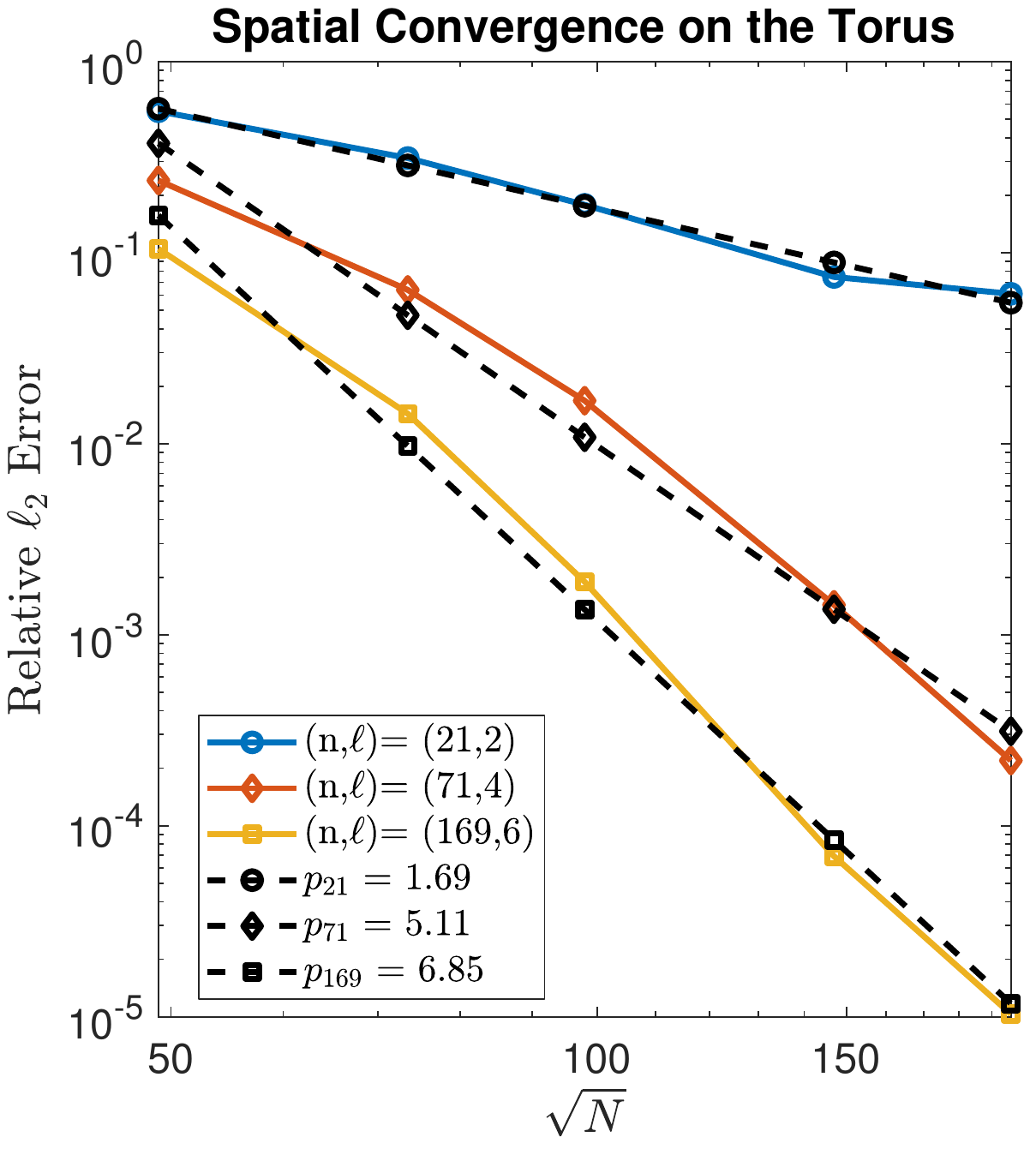}
\caption{Convergence on the torus for the surface advection equation of two cosine bells in a steady flow (left) and two Gaussian bells in the same flow (right). The figure shows relative $\ell_2$-error as a function of $\sqrt{N}$ for different values of stencil $n$ and polynomial degree $\ell$. The dashed lines indicate lines of best fit, and their slopes are indicated in the legend as a measure of convergence rates.}
\label{fig:toradv}
\end{figure}

\paragraph{Cosine bells}
The first test case is the advection of a pair of compactly supported cosine bells given by the following initial condition:
\begin{align*}
c(\vx,0) = 0.1 + 0.9 (q_1 + q_2),
\end{align*}
where $q_1$ and $q_2$ are compactly-supported functions given by:
\begin{align*}
q_{1,2} = 
\begin{cases}
\frac{1}{2}\lf(1 + \cos\lf(2\pi r_{1,2}(\vx) \rt)\rt) & \text{if}\ r_{1,2}(\vx) < 0.5, \\
0 & \text{if}\ r_{1,2}(\vx) \geq 0.5,
\end{cases}
\end{align*}
and $r_{1,2} (\vx) = \|\vx -\vp_{1,2}\|$.  As on the sphere, the solution has only one continuous derivative. 
The convergence results for advecting the cosine bells are shown in Figure \ref{fig:toradv} (left), which clearly shows that increasing the polynomial degree $\ell$ does not increase convergence rates when advecting the cosine bell, but again does improve the accuracy (as on the sphere).

\paragraph{Gaussian bells}
The second test case involves advecting a pair of Gaussian bells given by the following initial condition:
\begin{align*}
c(x,y,z,0) = e^{-a \lf(x + (1+1/3) \rt)^2 + y^2 - 1.5a z^2} + e^{-a \lf(x - (1+1/3) \rt)^2 + y^2 - 1.5a z^2},
\end{align*}
where $a=20$. This centers two Gaussians at the same locations as the cosine bells from the first test case, with the rapid falloff ensuring locality without destroying smoothness. 
The convergence results in the $\ell_2$-norm are shown in Figure \ref{fig:toradv} (right), which again shows that increasing $\ell$ gives close to predicted convergence rates when the solution is smooth.

\subsection{Surface Advection-Diffusion}
\begin{figure}[h!]
\centering
\begin{tabular}{cc}
\includegraphics[scale=0.65]{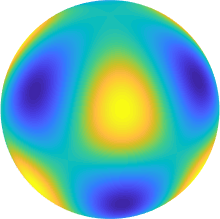}\hspace{1cm} & 
\includegraphics[scale=0.65]{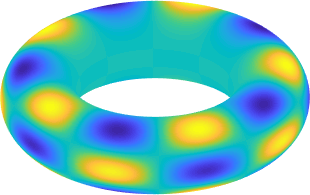} \\
(a) & (b)
\end{tabular}
\caption{Manufactured solutions at time $t=10^{-3}$ for surface advection-diffusion tests: (a) sphere \& (b) torus.}
\label{fig:initcond}
\end{figure}
Next, we investigate the convergence properties of our stabilized RBF-FD method on the surface advection-diffusion equation. Once again, we focus on the sphere and the torus, using the same node sets as in the pure advection case. In this set of tests, we simply check convergence against a manufactured solution on both surfaces. However, since we are now approximating both first and second-order differential operators, we set $\ell = \xi+1$, use $\xi=2,3,4$, and set all other parameters according to Algorithm \ref{alg:hyp}. All timestepping was done using the semi-implicit backward differentiation formula of order 4 (SBDF4)~\cite{Ascher97}, \revone{which is an IMEX method}. While this is nowhere near optimal for high Peclet numbers, we use it to illustrate that our hyperviscosity formulation can be stepped implicitly for efficiency. In this set of tests, we use Peclet numbers of $1$ and $100$ to study the stability and convergence of our method in different parameter regimes, mimicking similar tests for Euclidean domains from~\cite{SFJCP2018}. For these tests, we select the smaller time step between one that corresponds to a Courant number of 0.3, and one that attempts to keep the temporal error below the spatial error. In other words, the time step $\triangle t$ is given by:
\begin{align*}
\triangle t = \min \lf(0.3 N^{-\frac{1}{2}},  N^{-\frac{\xi}{8}}\rt),
\end{align*}
where $N$ is the number of nodes, $\xi$ is the desired spatial order of convergence, and the factor of $\xi/8$ comes from the fourth-order convergence of SBDF4 scheme. Here, we assume the nodes are quasi-uniform with a spacing of approximately $N^{-\frac{1}{2}}$.

\subsubsection{Advection-diffusion on the sphere}
\begin{figure}[h!]
\centering
\includegraphics[scale=0.61]{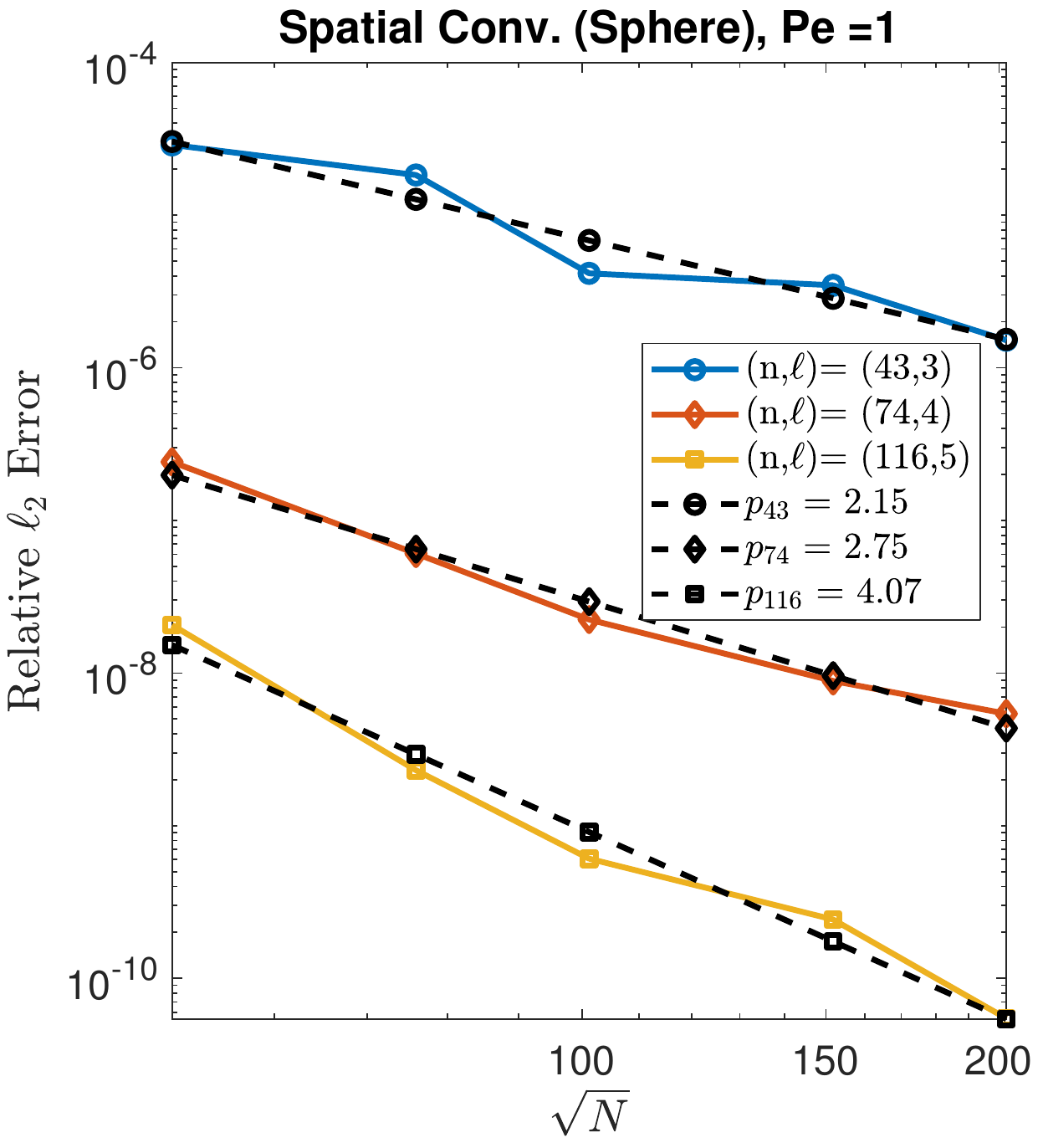}
\includegraphics[scale=0.62]{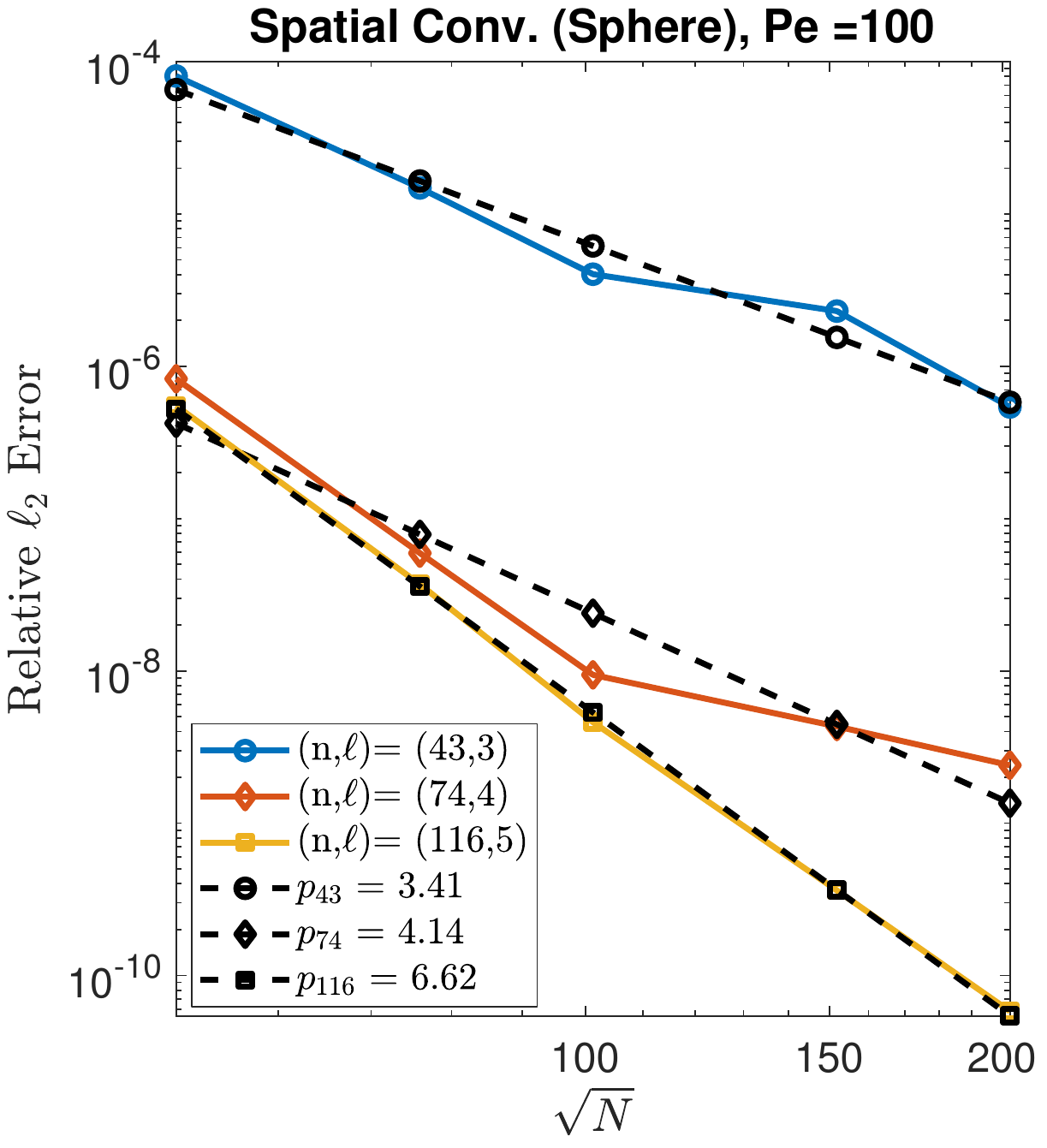}
\caption{Convergence on the sphere for the surface advection-diffusion equation at Peclet numbers $Pe=1$ (left) and $Pe=100$ (right). The figure shows relative $\ell_2$-error as a function of $\sqrt{N}$ for different values of stencil $n$ and polynomial degree $\ell$. The dashed lines indicate lines of best fit, and their slopes are indicated in the legend as a measure of convergence rates.}
\label{fig:sphadvdiff}
\end{figure}
For this test, we prescribe a solution and then manufacture a forcing term that makes the solution hold true. Our prescribed/manufactured solution is generated using the $Y_4^3$ real-valued spherical harmonic; more specifically, we shift it and multiply it by a time-dependent term. This is given by:
\begin{align*}
c(x,y,z,t) = 1 + \frac{3}{4}\sqrt{\frac{35}{2\pi}} \lf(x^2 - 3y^2\rt) x z \sin(t).
\end{align*}
The forcing term is computed analytically as:
\begin{align*}
f = \frac{\partial c}{\partial t} + \vu \cdot \nabla_{\mathbb{S}^2} c - \nu \Delta_{\mathbb{S}^2} c,
\end{align*}
where $\Delta_{\mathbb{S}^2} c = -20 (c - 1)$; the initial condition is shown at time $t=10^{-3}$ in Figure \ref{fig:initcond} (left). It is important to note that the forcing term does \emph{not} contain the hyperviscosity operator. This allows us to correctly verify the impact of the hyperviscosity formulation on convergence to the manufactured solution. The velocity field $\vu$ is the same as the steady velocity field from \eqref{eq:steady_vel}, with $\alpha = \pi/2$. We ran the simulation out to final time $T=2\pi$ (one full period). The results are shown in Figure \ref{fig:sphadvdiff}. Going from left to right across Figure \ref{fig:sphadvdiff}, we see that convergence rates increase as the Peclet number is increased. This occurs due to our choice of $\ell = \xi+1$. As we increase the Peclet number, the influence of the advection operator on transport increases relative to that of the diffusion operator. Since the advection operator is a first-order differential operator, this leads to an extra order of accuracy for the mixed character PDE since approximating the advection operator only requires a choice of $\ell = \xi$ (as seen in the previous subsection).

\subsubsection{Advection-diffusion on the torus}
\begin{figure}[h!]
\includegraphics[scale=0.62]{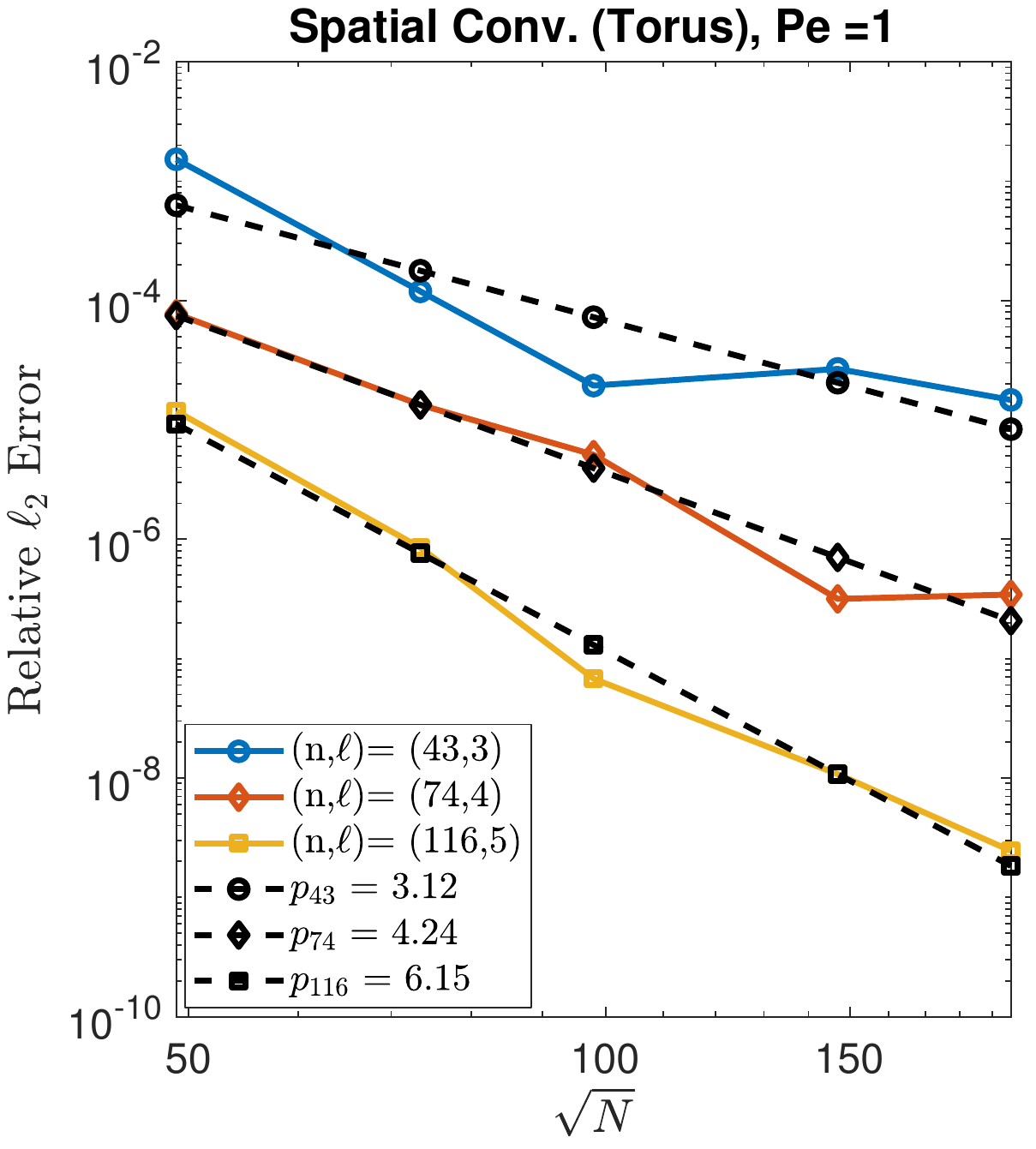}
\includegraphics[scale=0.62]{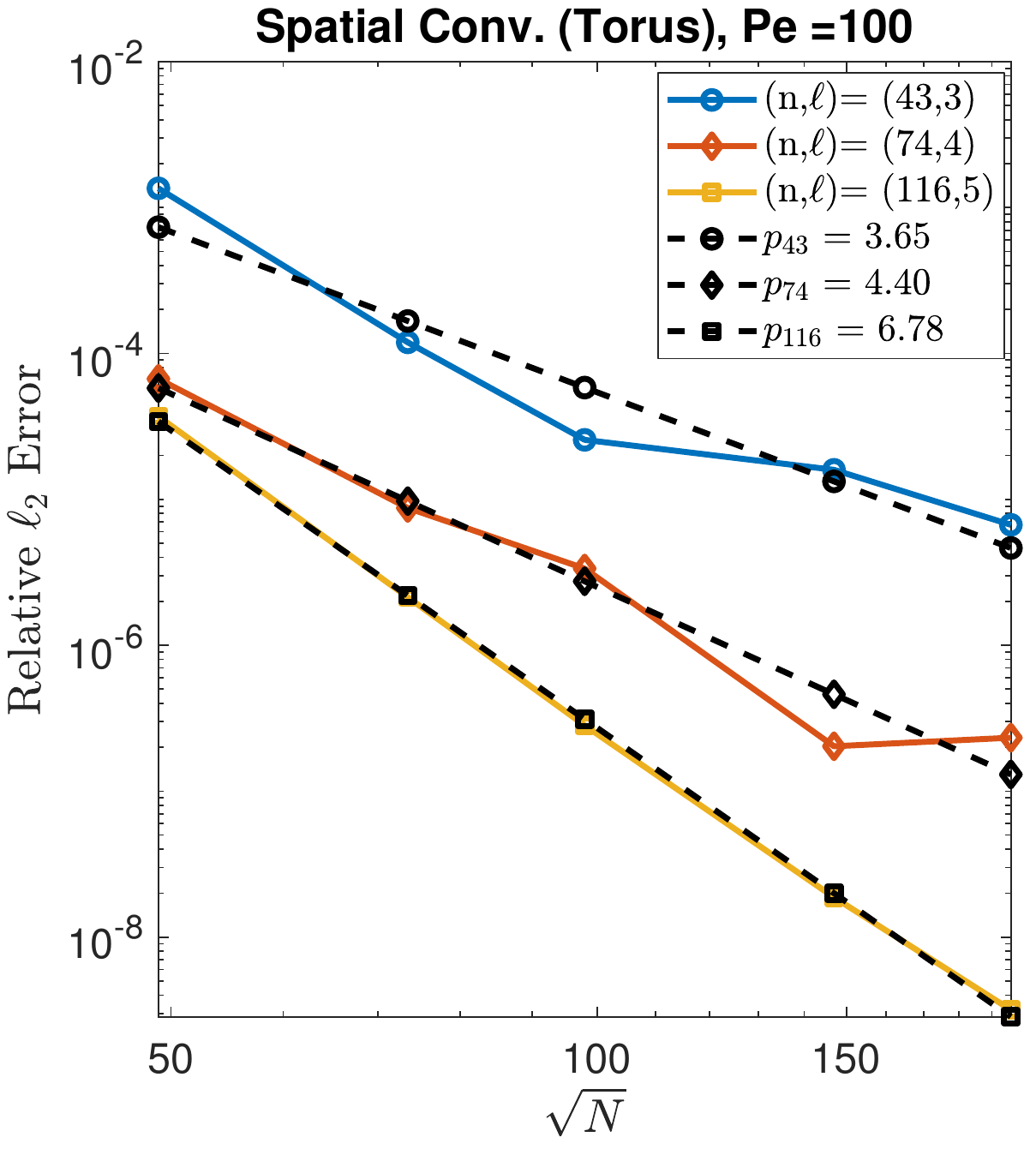}
\caption{Convergence on the torus for the surface advection-diffusion equation at Peclet numbers $Pe=1$ (left) and $Pe=100$ (right). The figure shows relative $\ell_2$-error as a function of $\sqrt{N}$ for different values of stencil $n$ and polynomial degree $\ell$. The dashed lines indicate lines of best fit, and their slopes are indicated in the legend as a measure of convergence rates.}
\label{fig:toradvdiff}
\end{figure}
We again prescribe a solution and then manufacture a forcing term that makes the solution hold true. Our prescribed/manufactured solution is given by the following smooth function:
\begin{align*}
c(x,y,z,t) = 1 + \frac{1}{8}x\lf(x^4 - 10x^2y^2 + 5y^4\rt)\lf(x^2 + y^2 - 60z^2\rt)\sin(t).
\end{align*}
This initial condition is shown in Figure \ref{fig:initcond} (right). The forcing term is computed analytically as:
\begin{align*}
f = \frac{\partial c}{\partial t} + \vu \cdot \nabla_{\mathbb{T}} c - \nu \Delta_{\mathbb{T}} c,
\end{align*}
where $\mathbb{T}$ represents the torus. Setting $\rho = \sqrt{x^2 + y^2}$, the surface Laplacian of the manufactured solution
\begin{align}
\Delta_{\mathbb{T}} c= \frac{-3}{8 \rho^2}x \lf(x^4 - 10x^2y^2 + 5y^4\rt)\lf(10248\rho^4 - 34335\rho^3 + 41359 \rho^2 - 21320\rho + 4000\rt) \sin(t).
\end{align}
Again, the forcing term does not account for the hyperviscosity operator. The velocity field $\vu$ is the same as the steady velocity field given by \eqref{eq:sv_tor1}. The simulation was run out to final time $T = \pi$. The results are shown in Figure \ref{fig:toradvdiff}. Going from left to right across Figure \ref{fig:toradvdiff}, we once again see that convergence rates increase as the Peclet number is increased, though they appear to be higher than predicted even at low Peclet number. This illustrates that our method is stable under spatial and order refinement, and that the hyperviscosity term vanishes appropriately thereby allowing us to recover the manfactured solution. Note that the implicit time-stepping of the hyperviscosity operator was done purely for illustrative purposes at the higher Peclet numbers. For sufficiently high Peclet numbers, it would likely suffice to use fully explicit time-stepping of all terms in the PDE with the RK4 method.

\section{Applications}
\label{sec:appl}

We now present applications of our hyperviscosity model to advection-diffusion-reaction systems on two surfaces: the red blood cell surface~\cite{SWFKJSC2014,FuselierWright:2013,SNKJCP2018} and Dupin's cyclide~\cite{FuselierWright:2013,SWFKJSC2014}. In all cases, we advance our simulations in time using the SBDF2 method~\cite{Ascher97} with a timestep of $\triangle t = 10^{-3}$.  For the spatial discretization, we use $\xi = 4$, and once again set all other parameters according to Algorithm \ref{alg:hyp}. \revone{Node sets on the cyclide were generated using the software Distmesh~\cite{distmesh}, and the algorithms from~\cite{SFKSISC2018} were used to generate node sets on the red blood cell. In the latter case, these consist of generalized spiral nodes on the sphere mapped to the red blood cell, with supersampling and decimation ensuring quasi-uniformity on the red blood cell.} 

\subsection{Advective Cahn-Hilliard on Dupin's cyclide}
\label{sec:cahn}
\begin{figure}[h!]
\centering
\begin{tabular}{cc}
\includegraphics[scale=0.5]{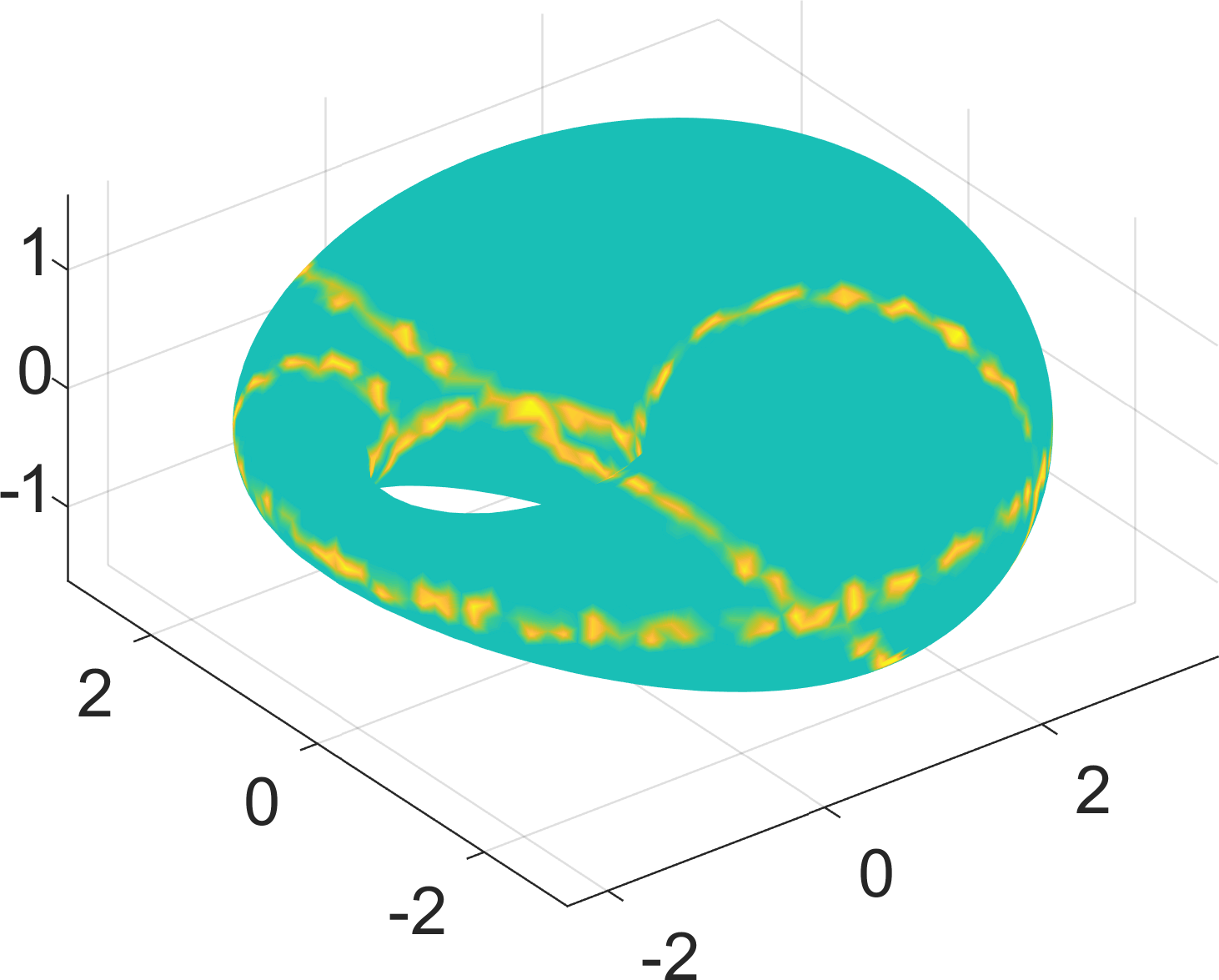} & 
\includegraphics[scale=0.5]{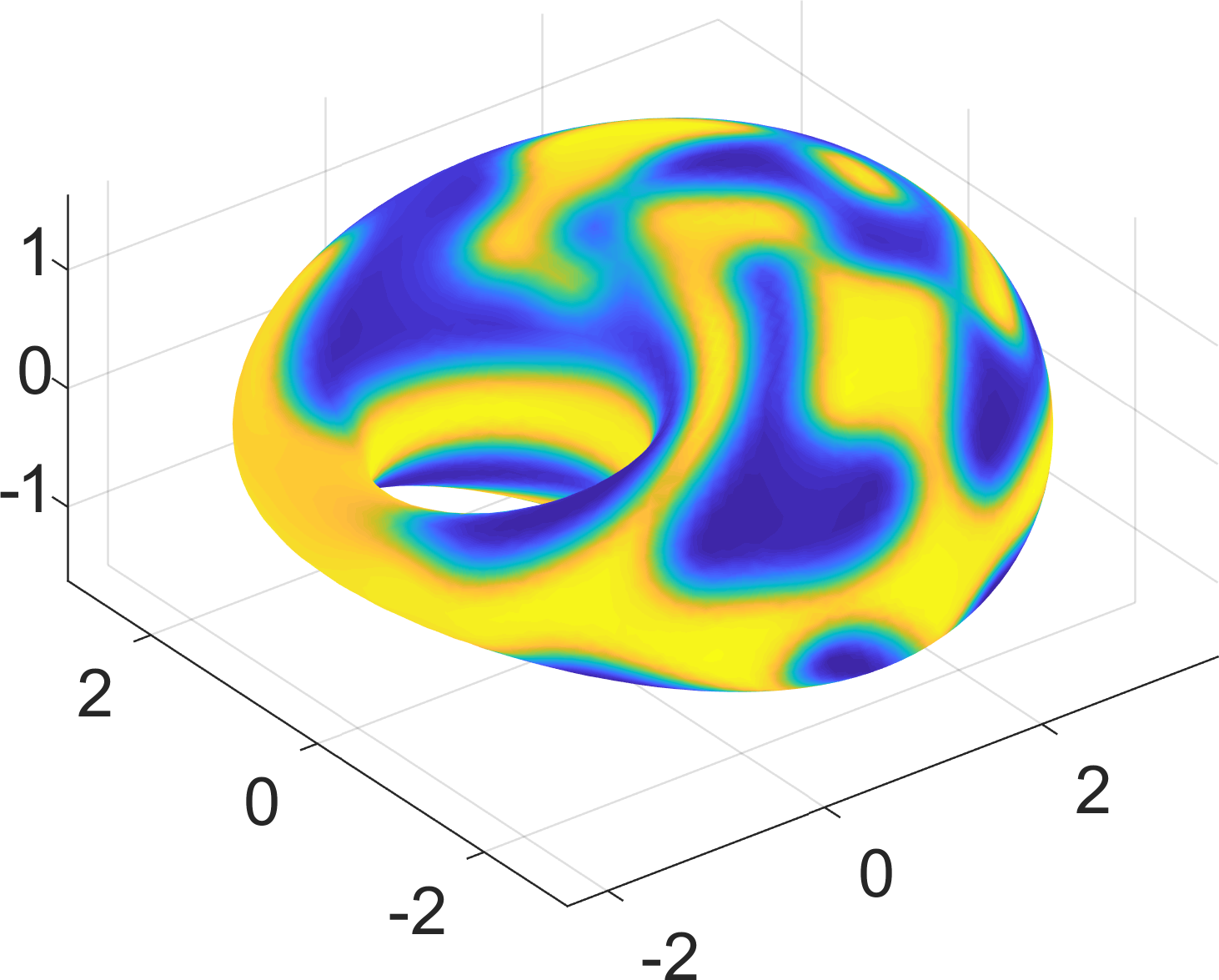} \\
\multicolumn{2}{c}{\includegraphics[scale=0.35]{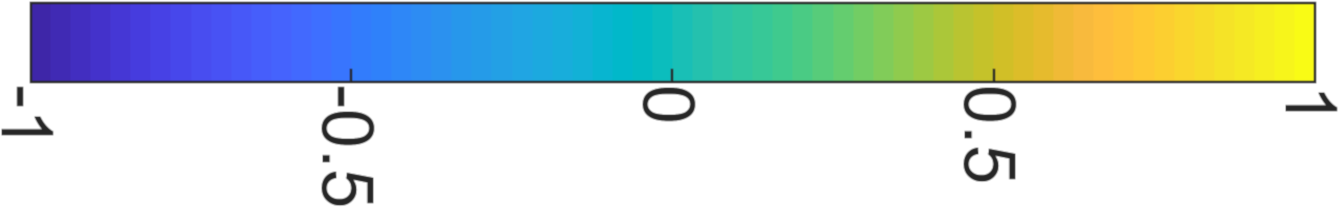}} \\
& \\
\includegraphics[scale=0.5]{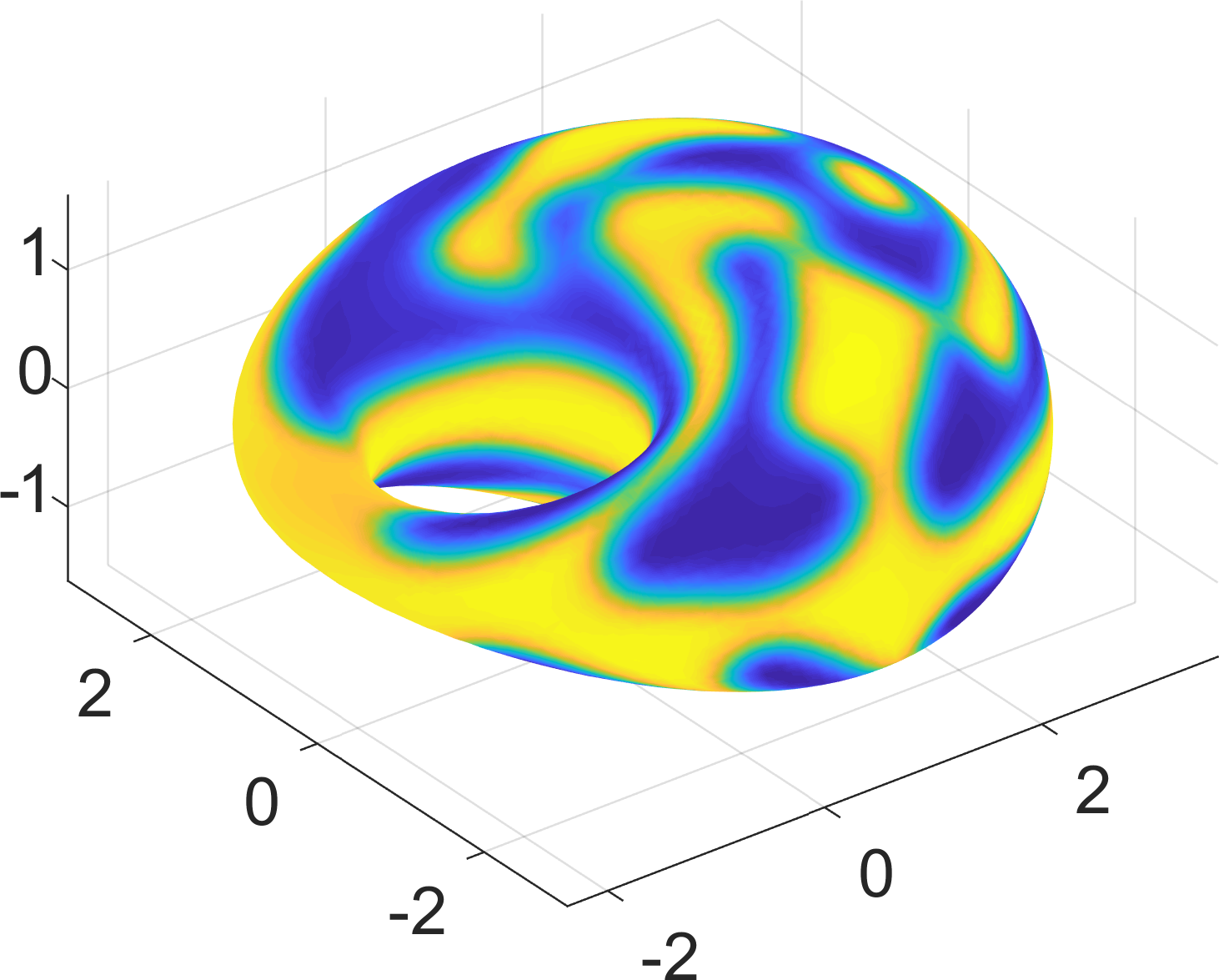} & 
\includegraphics[scale=0.5]{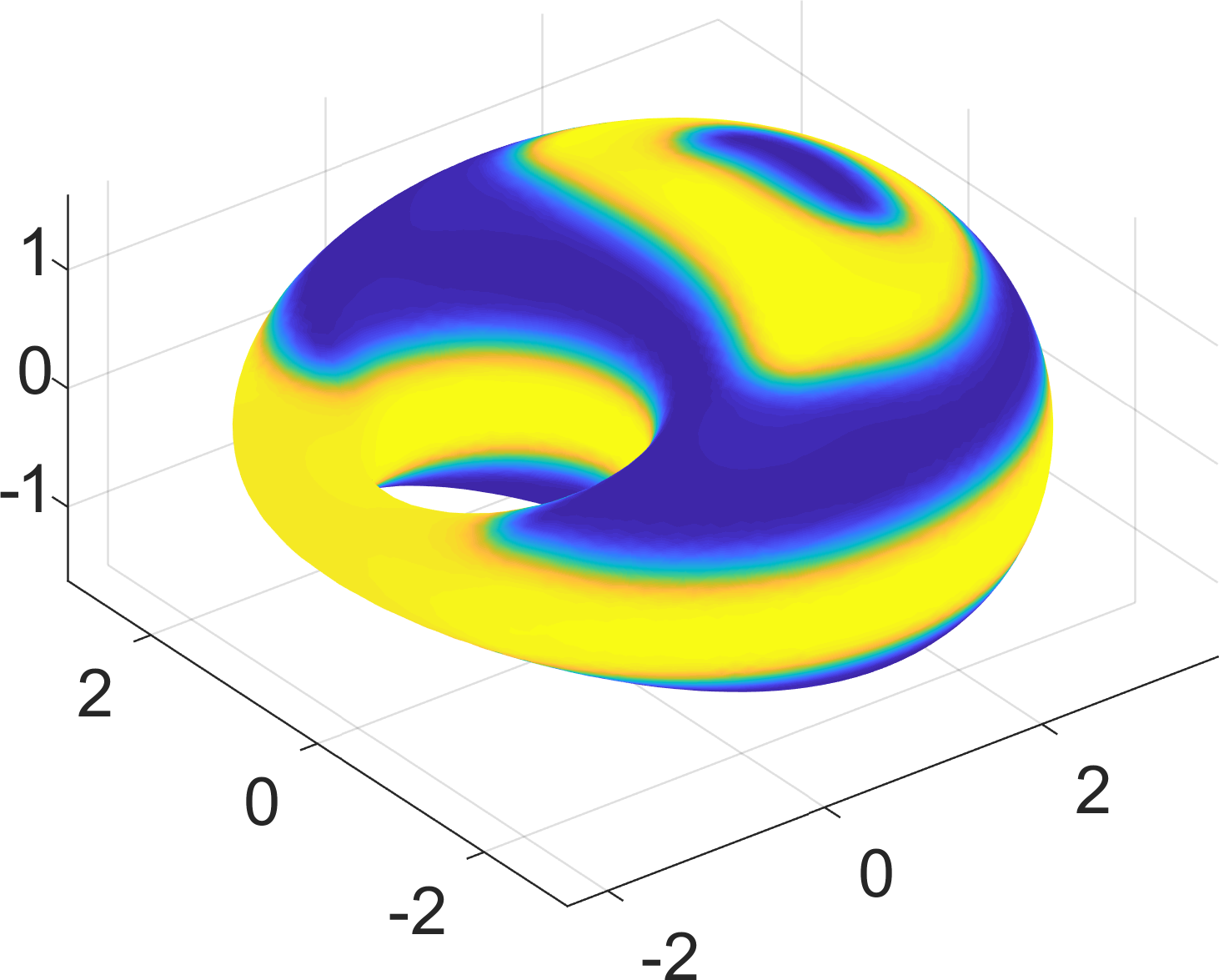} \\
& 
\end{tabular}
\caption{Solution of the advective Cahn-Hilliard equation \eqref{eq:chd} on Dupin's cyclide. The top left shows the initial condition, the bottom right the final solution at time $t=10$, and the other figures are intermediate snapshots.}
\label{fig:chd}
\end{figure}
We simulate an advective version of the Cahn-Hilliard equation on Dupin's cyclide. This is a nonlinear PDE governing phase separation with applications to both engineering and biology. The PDE (with artificial hyperviscosity) is given by:
\begin{align}
\frac{\partial c}{\partial t} + \nabla_{\mathbb{M}} \cdot \lf(c \vu \rt) = \nu \Delta_{\mathbb{M}} c^3 - \nu \Delta_{\mathbb{M}} c - \nu \beta \Delta^2_{\mathbb{M}} c + \gamma_1 \Delta^{\gamma_2}_{\mathbb{M}} c, \label{eq:chd}
\end{align}
where $\Delta^2_{\mathbb{M}}$ is the surface bilaplacian, and the $\gamma_1 \Delta^{\gamma_2}$ term is added solely to cancel the spurious growth in $ \nabla_{\mathbb{M}} \cdot \lf(c \vu \rt)$.  The velocity field $\vu$ is computed as the surface curl of a streamfunction $\psi$ whose gradient is given by:
\begin{align}
\nabla \psi &= -10 \cos\lf(\frac{\pi}{2} t\rt) \lf[1,0,0 \rt]^T, \label{eq:chvel1}\\
\vu &= \vn \times \nabla \psi. \label{eq:chvel2}
\end{align}
Since the velocity field $\vu$ is surface-incompressible ($\nabla_{\mathbb{M}} \cdot \vu = 0$) by construction, the solutions $c=1$ and $c=-1$ are critical points of the advection-diffusion-reaction system, causing the initial condition to separate over time into these two phases. We simulate \eqref{eq:chd} on Dupin's cyclide to time $t=10$ using a random initial condition (Figure \ref{fig:chd}a), with $\nu = 0.5$ and $\beta = 0.02$. As in~\cite{SNKJCP2018}, we approximate the surface bilaplacian as $B = L^2$, where $L$ is the differentiation matrix for the surface Laplacian. To balance stability and efficiency, we step the nonlinear terms $\nabla_{\mathbb{M}} \cdot \lf(c \vu \rt)$ and $\Delta_{\mathbb{M}} c^3$ explicitly in time, and step all other terms implicitly. The results for \revtwo{$N = 8266$} nodes are shown in Figure \ref{fig:chd}. While the advection term causes the phases to occasionally remix, they separate unless forced to remix by the advection term. However, we observed that the phases remained separated without remixing for a sufficiently large final time.
\begin{figure}[h!]
\centering
\includegraphics[scale=0.3]{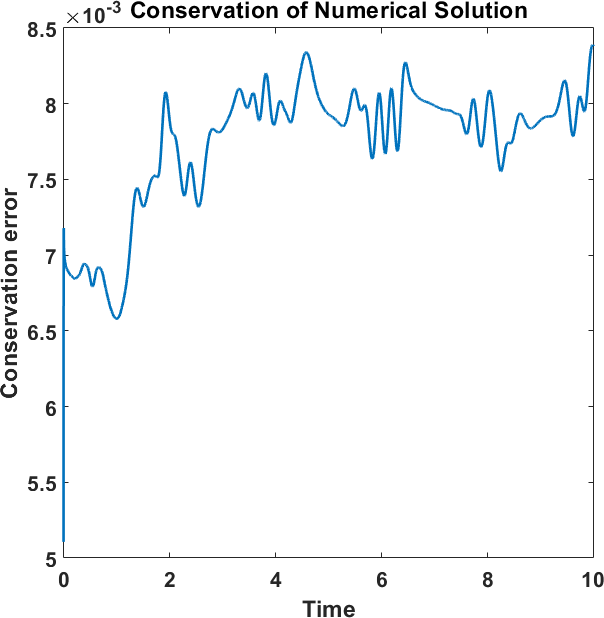}  
\caption{Conservation error in the numerical solution to the advective Cahn-Hilliard equation \eqref{eq:chd} on Dupin's Cyclide as a function of time.}
\label{fig:ch_cons}
\end{figure}
\revtwo{We also compute the conservation errors in the numerical solution as a function of time. If $c$ is the numerical solution and $c_0$ is its initial value at time $t=0$, the conservation error is given as $|\int_{\mathbb{M}} \lf(c - c_0 \rt)|$. The integral over the surface $\mathbb{M}$ is computed using an eighth-order accurate RBF-based quadrature technique from~\cite{ReegerFornbergQuad2}. The conservation error is shown as a function of time in Figure \ref{fig:ch_cons}. The figure shows that the conservation error grows in time. However, we have observed that increasing the order $\xi$ or the number of nodes $N$ significantly improves these errors, as seen in~\cite{SWJCP2018}. If necessary, the time-stepping can be modified to ensure conservation~\cite{Ling1}, but such an experiment is beyond the scope of this work.}

\subsection{Turing spots with advection on the red blood cell}
\label{sec:rbc} 
\begin{figure}[h!]
\centering
\begin{tabular}{cc}
\includegraphics[scale=0.5]{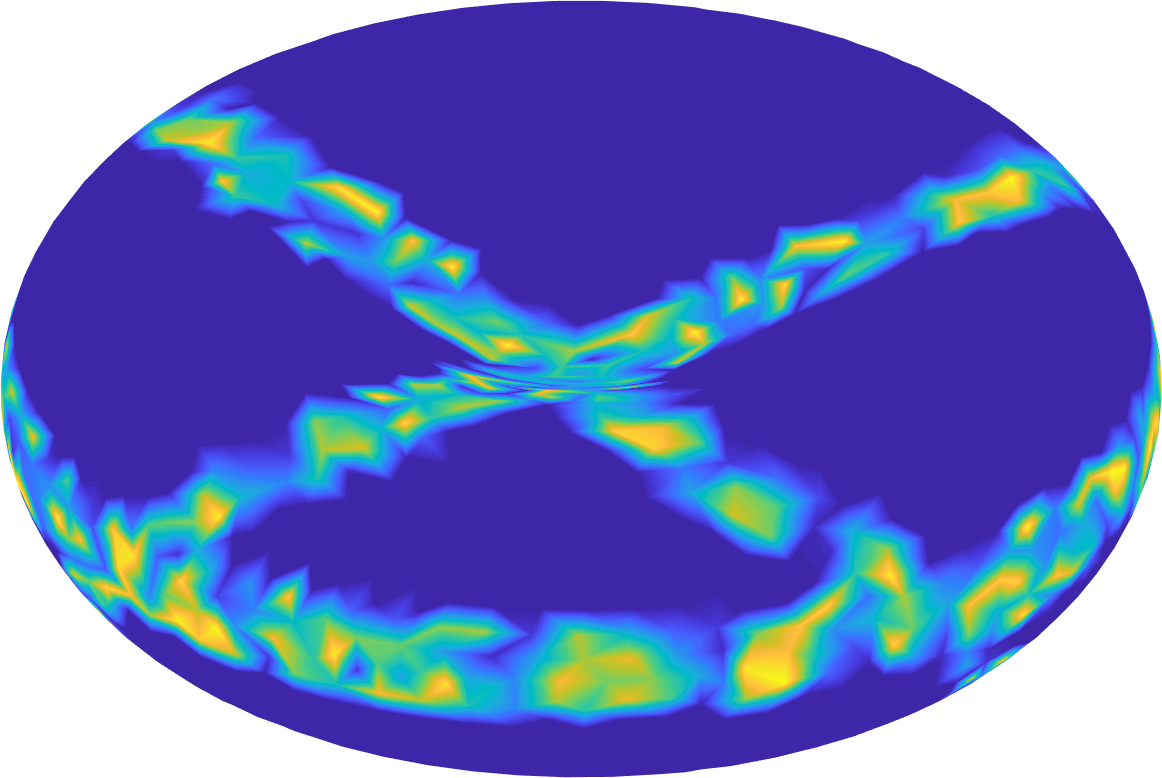} & 
\includegraphics[scale=0.5]{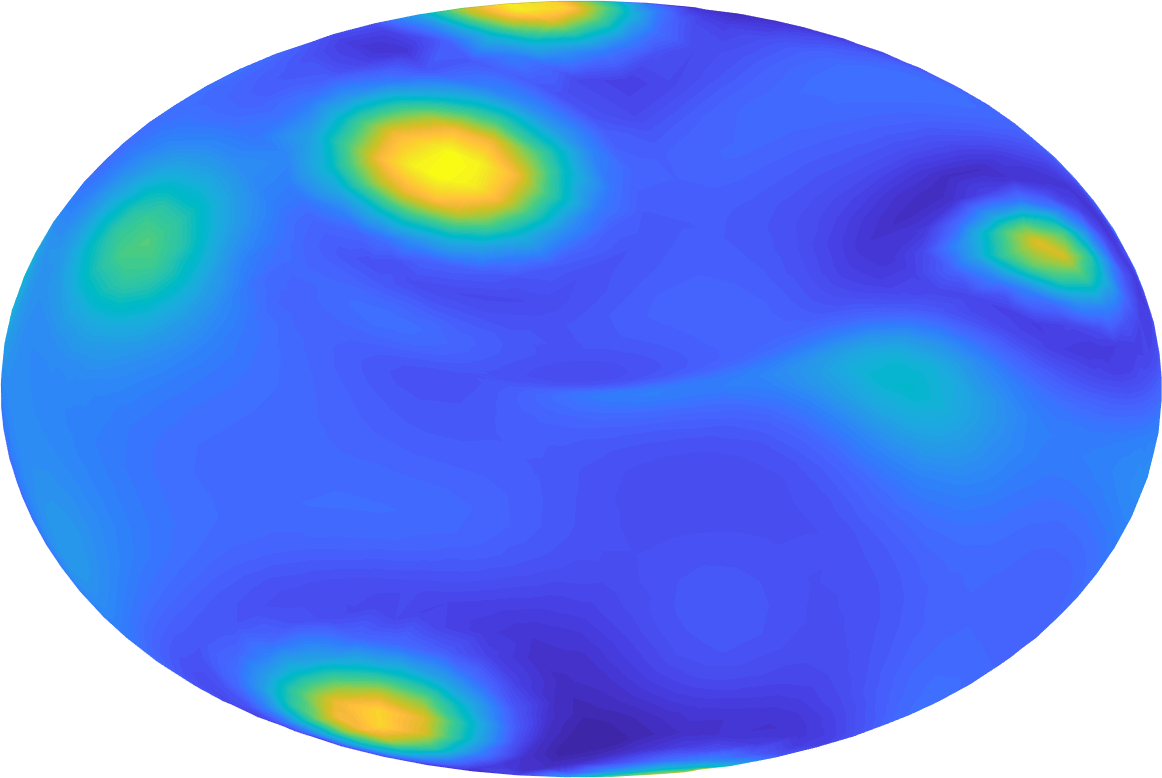} \\
& \\
\includegraphics[scale=0.5]{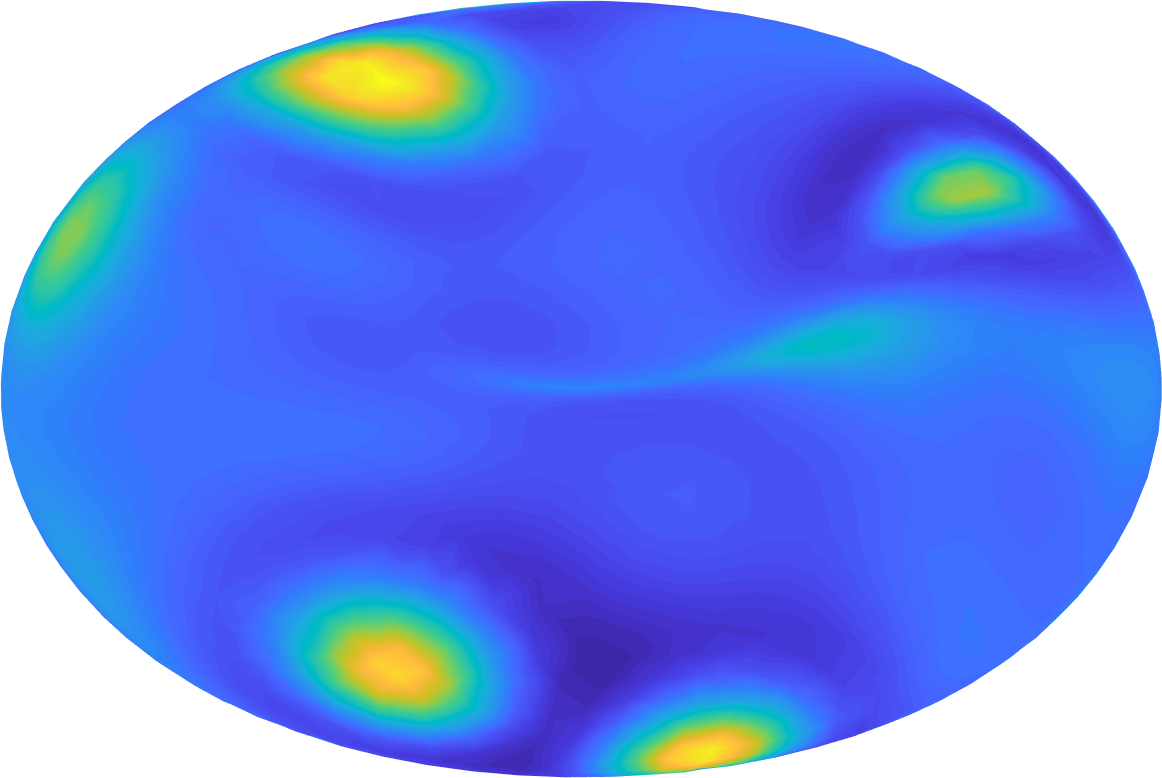} &
\includegraphics[scale=0.5]{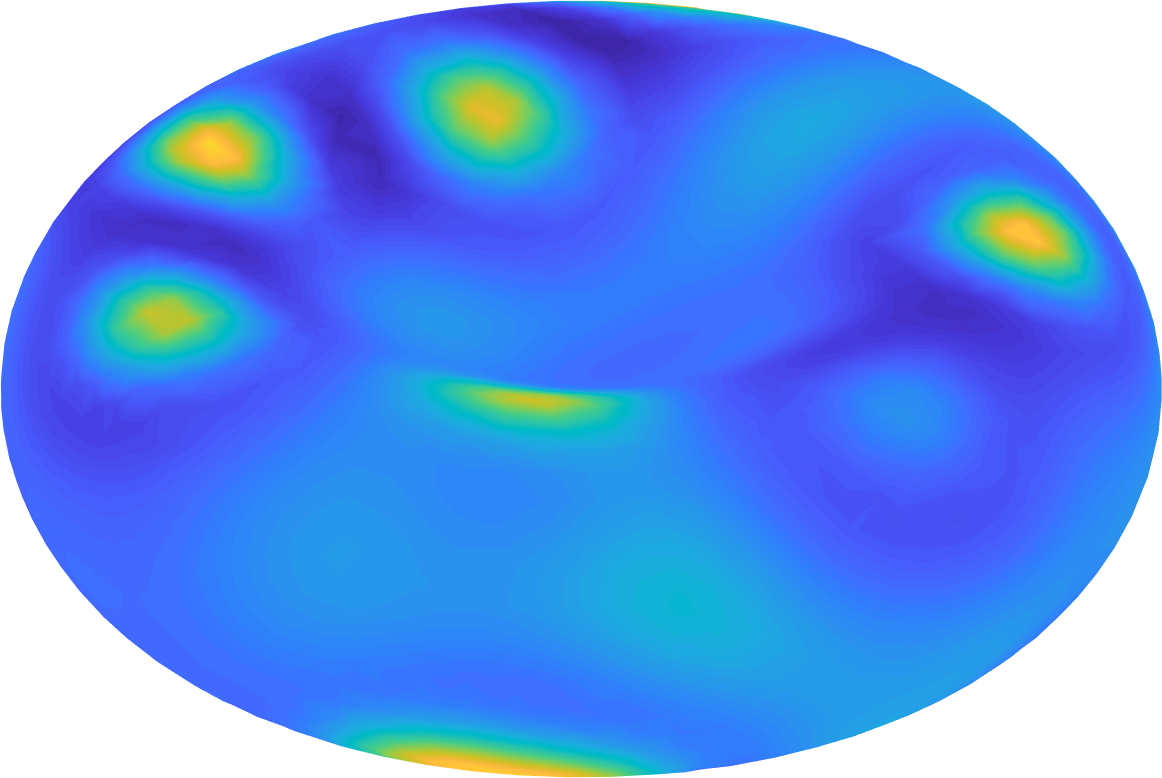} \\
& 
\end{tabular}
\caption{Solution of the advective Turing equations \eqref{eq:t1}--\eqref{eq:t2} on an idealized red blood cell. The top left shows the initial condition, the bottom right the final solution at time $t=800$, and the other figures are intermediate snapshots. A brighter color indicates a higher concentration of $c_1$, and a darker color indicates a lower concentration.}
\label{fig:rbc}
\end{figure}
Finally, we solve a (reaction-)coupled advection-diffusion-reaction system on an idealized red blood cell. Specifically, we simulate the pattern-generating Turing system given by:
\begin{align}
\frac{\partial c_1}{\partial t} + \nabla_{\mathbb{M}} \cdot \lf(c_1 {\bf w} \rt)&= \delta_1 \Delta_{\mathbb{M}}c_1 + \alpha c_1 \lf(1 - \tau_1 c_2^2\rt) + c_2 \lf(1-\tau_2 c_1\rt) + \gamma_1 \Delta^{\gamma_2}_{\mathbb{M}} c_1, \label{eq:t1} \\
\frac{\partial c_2}{\partial t} + \nabla_{\mathbb{M}} \cdot \lf(c_2 {\bf w} \rt)&= \delta_2 \Delta_{\mathbb{M}}c_2 + \beta c_2 \lf(1 + \frac{\alpha \tau_1}{\beta} c_1 c_2\rt) + c_1 \lf(\eta_1 + \tau_2c_2 \rt) + \gamma_1 \Delta^{\gamma_2}_{\mathbb{M}} c_2, \label{eq:t2}
\end{align}
where the hyperviscosity terms are only added to cancel spurious growth modes in the surface gradient terms. We use the parameters $\delta_1 = 0.0011$, $\delta_2 = 0.0021$, $\tau_1 = 0.02$, $\tau_2 = 0.2$, $\alpha = 0.899$, $\beta = -0.91$, and $\eta_1 = -\alpha$, with the initial condition as shown in Figure \ref{fig:rbc}a. This set of parameter choices promotes spot formation in the absence of advection~\cite{SWFKJSC2014,SNKJCP2018}. The velocity field ${\bf w}$ is simply a scaled version of \eqref{eq:chvel1}--\eqref{eq:chvel2}, \emph{i.e.}, we set ${\bf w} = 0.1 \vu$. We chose this scaling because the Turing patterns (specifically, spots) are sensitive to the choice of the velocity field, with a very fast velocity field inhibiting spot formation. As before, the advection and reaction terms were stepped explicitly in time, while the hyperviscosity and diffusion terms are stepped implicitly. The simulations were run out to a final time of $t= 800$, and the results \revone{$N=2553$ nodes} are shown in Figure \ref{fig:rbc}. We see both the formation of spots from the initial condition in Figure \ref{fig:rbc}b, and their motion due to the velocity field (going from Figure \ref{fig:rbc}b to Figure \ref{fig:rbc}c).

\section{Summary and Future Work}
\label{sec:summary}

We have presented a novel, automatic procedure for stabilizing discretizations of linear advection and advection-diffusion equations on manifolds. In particular, we have investigated the procedure in the context of overlapping RBF-LOI discretizations. The geometric, high-order flexibility of RBF-LOI discretizations when paired with automatically-parameterized hyperviscosity yields a numerical algorithm that is efficient, accurate, and locally and globally stable for all examples we have investigated. We have tested our algorithm on surface advection and surface advection-diffusion problems on manifolds of co-dimension-1 embedded in $\mathbb{R}^3$. The automated choice of hyperviscosity essentially eliminates tedious and expensive hand-tuning of stabilization parameters when applying our RBF-LOI algorithm across variegated PDEs and geometries. 

The mathematical strategy we have developed for automated tuning of the hyperviscosity term provides general guidance for stabilizing both linear and linearized PDEs on manifolds, can be easily applied to non-RBF-type discretizations, and is easily generalized to manifolds embedded in arbitrary Euclidean dimensions. Our analysis also provides diagnostic guidance for determining when discretized advection-diffusion systems possess sufficient diffusion to obviate the need for stabilizing hyperviscosity terms.

\section*{Acknowledgments}

VS was supported by NSF grants DMS-1521748 and CISE 1714844. GBW acknowledges funding support for this project under NSF CCF 1717556. AN was partially supported by NSF DMS-1720416 and AFOSR FA9550-15-1-0467. 

%\appendix
%\input{Appendix}

%\section*{References}
\bibliographystyle{siamplain}
\bibliography{article_refs_mod}

\end{document}